\documentclass{amsart}
\usepackage{url}
\usepackage{amssymb}
\usepackage{mathtools}
\usepackage{hyperref}

\newtheorem{theorem}[equation]{Theorem}
\newtheorem{prop}[equation]{Proposition}
\newtheorem{lemma}[equation]{Lemma}
\newtheorem{corollary}[equation]{Corollary}
\theoremstyle{definition}

\newtheorem{remark}[equation]{Remark}
\newtheorem{example}[equation]{Example}
\theoremstyle{plain}
\numberwithin{equation}{section}
\numberwithin{figure}{section}

\DeclareMathOperator{\im}{Im}
\DeclareMathOperator{\re}{Re}
\DeclareMathOperator{\Gal}{Gal}
\DeclarePairedDelimiter\abs{\lvert}{\rvert}

\def\mybullet{\noindent$ \bullet $ }

\let\myv\v % for Nekovar's name

\def\a{\alpha}
\def\b{\beta}
\def\g{\gamma}
\def\l{\lambda}
\def\v{\varepsilon}
\def\s{\sigma}
\def\t{\tau}

\def\C{\mathbb C}
\def\H{\mathcal H}
\def\P{\mathbb P}
\def\Q{\mathbb Q}
\def\R{\mathbb R}
\def\Z{\mathbb Z}

\def\DD{\mathcal{D}}
\def\OO{\mathcal{O}}

\def\cA{\mathcal{A}}
\def\cB{\mathcal{B}}
\def\cC{\mathcal{C}}
\def\cE{\mathcal{E}}
\def\cF{\mathcal{F}}

\def\PP{\mathfrak{P}}
\def\tQ{\widetilde Q}

\def\W{W} % for the remark at the end of Section 3.
\def\Wint{\W_{\textup{int}}}

\def\ge{\geqslant}
\def\le{\leqslant}

\def\mycoprod{\raise1pt\hbox{\,$ \scriptstyle\coprod $\,}}

\def\INTE{K_2^T(E)_{\textup{int}}}
\def\INTC{K_2^T(C)_{\textup{int}}}

\def\vph{\vphantom{$b^{b^{b^b}}p_{p_{p_{p_p}}}$}} % for some vertical spacing in tables.

\begin{document}

\title[$K_2$ of elliptic curves over non-Abelian fields]{$\textit{K}_{\mathbf 2}$ of families of elliptic curves over non-Abelian cubic and quartic fields}

\author[F. Brunault]{Fran\c{c}ois Brunault}
\address{ENS Lyon, UMPA, 46 all\'{e}e d'Italie, 69007 Lyon, France}
\email{francois.brunault@ens-lyon.fr}

\author[R. de Jeu]{Rob de Jeu}
\address{Faculteit der B\`etawetenschappen\\Afdeling Wiskunde\\Vrije Universiteit Amsterdam\\De Boelelaan 1111\\1081 HV Amsterdam\\The Netherlands}
\email{r.m.h.de.jeu@vu.nl}

\author[H. Liu]{Hang LIU$^{\ast}$}
\address{School of Mathematical Sciences, Shenzhen University,
Shenzhen, 518060, Guangdong, P. R. China}
\email{liuhang@szu.edu.cn}
\address{$^{\ast}$Corresponding author}

\author[F. Rodriguez Villegas]{Fernando Rodriguez Villegas}
\address{The Abdus Salam International Centre for Theoretical Physics, Stada Costiera 11,
Trieste 34151, Italy}
\email{villegas@ictp.it}

\begin{abstract}
We give two constructions of families of elliptic curves over cubic or
quartic fields with three, respectively four, `integral' elements in the
kernel of the tame symbol on the curves.
The fields are in general non-Abelian, and the elements linearly independent.
For their integrality, we prove a new criterion that does not ignore any torsion.
We also verify Beilinson's conjecture numerically for just over 90 of
the curves.
\end{abstract}

\keywords{$K_2$, integrality, regulator, Beilinson's conjecture, elliptic curve}
\subjclass[2010]{Primary: 19E08, 19F27; secondary: 11G40}

\maketitle

\section{Introduction}

Based on pioneering results by Borel~\cite{Borel77} for number
fields, and Bloch~\cite{Bl00} for certain elliptic curves over~$ \Q $,
Beilinson in~\cite{Be85} stated very far-reaching
conjectures on the relation between special values of $L$-functions
and regulators of $K$-groups of smooth projective varieties
defined over number fields.
These conjectures give precise information on these
$K$-groups (e.g., the maximum number of linearly independent
elements in such a $ K $-group) that is currently out of reach in general.

To verify the conjecture (theoretically or numerically) one usually
constructs as many linearly independent elements in the relevant
$K$-group as the conjecture predicts, computes their images under the
Beilinson regulator map and relates those to the corresponding $L$-value.

On the theoretical side, this has been carried out by Beilinson
in the case of modular curves \cite{Be86}. Others have considered
different situations (see, e.g.,~\cite{Ots, Lem17, CLR, KLZ, BC16}),
where the Beilinson regulator is usually expressed as an automorphic integral
and can sometimes be related to the relevant $L$-value.
We refer to the survey papers~\cite{DS91, Nek94}
for a detailed exposition of Beilinson's conjecture
and for older (and by now classical) results along those lines.
Numerically, checks of the Beilinson conjecture for curves were carried out in, e.g.,
\cite{Young,dJ96,DJZ,Asa18}.

In this article, we study elliptic curves $E$ over certain number fields $F$
of fixed degree, themselves varying in a family.
In our setting the relevant $K$-group, which we denote by $\INTE$,
is the subgroup of $K_2$ of the function field of $E$
consisting of elements having trivial tame symbol at every point of $E$,
and satisfying a certain integrality condition; see~\cite[\S1]{LJ15}
for the necessary definitions and~\cite[\S3]{DJZ} for more background.
Beilinson's conjecture predicts that the $\Z$-rank\footnote{This
formulation actually supplements Beilinson's conjecture with
a conjecture by Bass on finite generation; see the footnote on~\cite[p.343]{DJZ}.}
of this $K$-group
is equal to the degree of $F$; this is currently not known for a single
elliptic curve. In fact, we are not aware of any examples in the literature
of elliptic curves~$ E $ over non-Abelian~$ F $ for which it is known
that~$ \INTE $ has at least that many linearly independent elements.

The main focus in this paper is to give two constructions of families of elliptic curves
over (mostly non-Abelian) cubic or quartic fields with 3 and 4 linearly independent
elements in $\INTE$ respectively, thus matching the ranks
predicted by the conjecture in an entirely new situation.

We first construct elements in the kernel $K_2^T(E)$ of the tame symbol using a method of Bloch,
taking as input functions with divisors in a finite subgroup of $E$
(see Section~\ref{blochsection}).
Doing this for the universal families~$E_t$ of elliptic curves having a point of order $N=7$, $8$ or $10$,
we find potentially linearly independent elements $S_1, \ldots, S_d$ in $K_2^T(E_t)$,
with~$d = 3$, $3$ or~$4$ respectively, and where we can specialise
the parameter~$ t $ to values in number fields.
In order to investigate the integrality of the resulting elements
we present in Section~\ref{new-integrality} a new
criterion (not ignoring any torsion) for the integrality of certain
elements of $K_2^T(C)$ for a curve~$ C $
over a number field (see Corollary~\ref{new-int-cor}).
Applying it to the above universal families, we find
that the elements $S_i$ are integral when the parameter~$t$
satisfies certain conditions akin to the unit equation
in number fields.
These conditions are satisfied for some explicitly described~$t$
that belong to certain degree~$d$ number fields, themselves parametrised
by an integer~$a$.

To establish the linear independence of the elements
$S_1, \ldots, S_d$, we compute the limit behaviour of the
Beilinson regulator of these elements when~$\abs{a}$
goes to infinity. In the case of the universal elliptic curves
above, this amounts to understand the behaviour of the regulator
near the cusps of the modular curve.

Let us here state a sample of our main result for~$N=7$, where~$ d = 3 $.

\begin{theorem} \label{intro-thm}
Let $ f_a(X) $ be one the following polynomials
\begin{align}
\label{fa 1}
& X^3 + a X^2 - (a+3) X + 1 \qquad (a \in \Z, \; a \neq -8), \\
\label{fa 2}
& X^3 + a X^2 - (a+1) X + 1 \qquad (a \in \Z).
\end{align}
Let $F = \Q(t)$ be the cubic field generated by a root $t$ of $f_a$.
Consider the elliptic curve~$E$ over~$F$ defined by
 \begin{equation*}
 E : y^2 + (1-g) xy - f y = x^3 - f x^2
 \end{equation*}
with $f = t^3-t^2$ and $g = t^2-t$.
Then~$ S_1 $, $ S_2 $ and~$ S_3 $ are in~$\INTE$,
and they are $\Z$-linearly independent for~$|a| \gg 0$.
\end{theorem}
The family \eqref{fa 1} is the family of \emph{simplest
cubic fields} introduced by Shanks~\cite{Sh}.
By contrast, the family \eqref{fa 2} defines non-Abelian
cubic fields.

The results for $ N = 8 $, where~$ d = 3 $, are similar. For~$ N = 10 $,
where~$ d= 4 $, we find many families of the corresponding quartic number fields, most
of which are non-Abelian. For the precise statements
we refer the reader to Theorem \ref{thm:independence}.

The other construction, using methods that are independent of the rest of the paper,
is given in Section \ref{K2-cubic}, where we consider other families of elliptic
curves $E$ defined over the same cubic fields as in the first construction.
Using non-torsion points on the elliptic curves,
we construct 3 elements in $\INTE$ that are linearly independent when the
parameter is large enough.

The second part of Beilinson's conjecture predicts that the
Beilinson regulator of the elements in~$ \INTE $ is a non-zero rational multiple of
the leading coefficient~$ L^*(E,0) $ in the Taylor expansion of $L(E,s)$ at $s=0$.
For elliptic curves over general number fields, there are very few
results in this direction.
Although we are not able to link the regulator with~$L^*(E,0)$ theoretically,
we check this numerically in a number of cases for both constructions.

Our results are of a different nature than Beilinson's theorem \cite{Be86},
which relates regulators for modular curves to $L$-values of modular forms.
For elliptic curves $E$ over $\Q$, this theorem
leads to inexplicit elements in~$ \INTE $ (see \cite[Theorem~7.3.1]{scsc88}),
whose regulators are related to $L^*(E,0)$.
This can be extended to strongly modular elliptic curves defined over Abelian fields~\cite{Br18}.
By contrast, our elements in~$ \INTE $ are explicit, and our elliptic curves
are in general defined over non-Abelian fields.

The structure of this paper is as follows.
In Section~\ref{conjecture} we recall Beilinson's conjecture
in the form that we shall use in this paper.
In Section~\ref{new-integrality} we prove our new integrality
criterion. This section is independent of
the rest of the paper, and its results apply to more
general curves and more general bases than rings of algebraic integers.
In Section~\ref{K2-families} we give our construction of elements
based on a point of order~$ N $ on the curve. We then discuss how our
new integrality condition as formulated in Corollary~\ref{new-int-cor}(2)
leads us, for~$ N = 7 $, 8 or~10, to consider specific cubic and quartic fields.
These fields are given explicitly in Lemma~\ref{lemma:unit}.
We then formulate our main result on this construction in Theorem~\ref{thm:independence}.
The linear independence part in this theorem is achieved by
computing the limit behaviour of the Beilinson regulator in our
families. As this is more involved, it is given in Section~\ref{section:independence}.
In Section~\ref{K2-cubic} we give a much simpler construction
of families of elliptic curves with three integral elements,
over the same cubic families as in Theorem~\ref{intro-thm},
and prove their linear independence (except for finitely many
curves), again by computing a limit result for their Beilinson
regulator. As mentioned above, this section is almost independent of
the rest of the paper.
Finally, in Section~\ref{numerics} we verify Beilinson's conjecture numerically for 93 of the curves
in Theorem~\ref{thm:independence} and Section~\ref{K2-cubic},
sometimes finding large prime factors in the numerator of the resulting rational number.

\section{Beilinson's conjecture} \label{conjecture}

In this section, we make precise the relation that Beilinson's conjecture
predicts between the~$ L $-function~$ L(E,s) $ of an elliptic curve $ E $ over a number field~$ F $,
and~$ \INTE $.  We shall follow~\cite[Remark~3.14]{DJZ}, which also deals with more general curves.

First, let $ E $ be an elliptic curve over~$ \C $, and let~$\a = \sum_j \{f_j, g_j \} $ be
in $K_2^T(E)$.
For~$\g$ in~$ H_1(E(\C),\Z)$, the \emph{regulator pairing}
between $ \g $ and $ \a $ is (well-)defined by
\begin{equation} \label{regulator-pairing}
\langle \g, \a \rangle = \frac{1}{2\pi}\int_{\g} \sum_j \eta(f_j, g_j)
\end{equation}
with $ \eta(f, g)=\log |f| \mathrm{d} \arg (g)-\log |g| \mathrm{d} \arg (f) $
for non-zero functions $ f $ and $ g $ on~$ E $, where we use
a representative of $ \g $ that avoids all zeroes and poles of
the functions involved.

Now let $ E $ be an elliptic curve defined over a number field~$ F $
of degree~$ m $.
Let $ X $ be the Riemann surface consisting of all~$ \C $-valued
points of $ E $, using all embeddings of $ F $ into~$ \C $.
It is a disjoint union of the complex points of~$ m $ elliptic curves~$ E^\s $,
obtained by applying each embedding~$ \s $ of~$ F $ into~$ \C $
to the coefficients of~$ E $.
Complex conjugation acts on $ X $ through its action on $ \C $.
If~$ \s $ is a real embedding this action is the usual action
on~$ E^\s(\C) $, with~$ E^\s $ defined over~$ \R $.
If~$ \s $ and $ \bar\s $ are two conjugate non-real embeddings of~$ F $,
then the action will interchange~$ E^\s(\C) $ and~$ E^{\bar\s}(\C) $.
From this one sees that $ H_1(X,\Z)^- $, the subgroup of~$ H_1(X,\Z) $
on which complex conjugation acts as multiplication by~$ -1 $,
is free of rank~$ m $. In fact, a $ \Z $-basis of it can be found
by combining $ \Z $-bases of~$ H_1(E^\s(\C), \Z)^- $
for all real embeddings~$ \s $ with $ \Z $-bases of
\begin{equation*}
H_1(E^\s(\C) \mycoprod E^{\bar\s}(\C) , \Z)^- = \left(  H_1(E^\s(\C) , \Z ) \oplus H_1( E^{\bar\s}(\C) , \Z) \right)^- \simeq  H_1(E^\s(\C) , \Z )
\end{equation*}
for all pairs $ \s, \bar\s $ of conjugate non-real embeddings.

We define a pairing
\begin{alignat*}{1}
H_1(X, \Z) \times K_2^T(E) & \to \R
\\
(\g , \a) & \mapsto \langle \g, \a \rangle_X
\end{alignat*}
as follows.
One has~$ \g = (\g_\s)_\s $ in~$ H_1(X, \Z) = \oplus_\s H_1(E^\s(\C) , \Z) $.
For each embedding~$ \s $ we pull back~$ \a $ to~$ E^\s $,
obtaining an element $ \a^\s $ in~$ K_2^T(E^\s) $. Note this
means~$ \a^\s $ is obtained from~$ \a $ by applying~$ \s $ to the coefficients of the functions in~$ \a $.
Then~$ \langle \g , \a \rangle_X = \sum_\s \langle \g_\s, \a^\s \rangle $ with
each pairing in the sum defined by~\eqref{regulator-pairing}.

Now let
\begin{equation*}
 \Lambda(E,s) = A^{s/2} (2 \pi)^{-sm} \Gamma(s)^m L(E,s)
\end{equation*}
with $ \Gamma(s) $ the Gamma-function and
$ A = N_{F/\Q}(f_E) d_{F}^2 $,
where~$ f_E $ is the conductor of~$ E/F $, and $ d_F $
the discriminant of~$ F $.
We assume~$ \Lambda(E,s) $ can be analytically continued
to the complex plane and satisfies a functional equation
$ \Lambda(E,2-s) = w \Lambda(E,s) $ with $ w = \pm 1 $, as stated in the Hasse-Weil conjecture
(see \cite[\S16.3]{husemoller}).
Then~$ L(E,s) $ has a zero of order~$ m $
at $ s = 0 $, and we let~$ L^*(E,0) = (m!)^{-1} L^{(m)}(E,0) $
be the first non-vanishing coefficient in its Taylor expansion
in~$ s $ at~$ 0 $.

Let $ \g_1,\dots,\g_m $ and~$ \a_1,\dots, \a_m $ form $ \Z $-bases of~$ H_1(X, \Z)^- $ and~$ \INTE $ modulo torsion
respectively, where the latter is supposed to be free of rank~$ m $
by the part of the conjectures mentioned in the introduction.
Then Beilinson expects that
\begin{equation} \label{eqn:L-reg}
L^*(E,0) = Q \cdot |\! \det ( \langle \g_i, \a_j \rangle_X  )_{i,j} |
\end{equation}
for some non-zero rational number~$ Q $. 
We call $R = |\! \det ( \langle \g_i, \a_j \rangle_X  )_{i,j} |$ the \emph{Beilinson regulator} of the elements $\a_j$.
In particular, this
implies that the pairing~$ \langle \,\, , \,\, \rangle_X $
between~$ H_1(X, \Z)^- $, and~$\INTE$  modulo torsion,
is non-degenerate.
The absolute value of the determinant is independent of the $ \Z $-bases used, making~$ Q $ well-defined
(with sign matching that of $ L^*(E,0) $).

\section{A new integrality criterion} \label{new-integrality}
For the elements in $ K_2^T(E) $ that we shall construct in Section~\ref{blochsection}
for certain elliptic curves~$ E $,
based on a rational point~$ P $ on~$ E $ of finite order,
we want to have a description of the `integrality condition'
that does not ignore any torsion. We achieve this in the case
that will be of interest to us, where,
on the minimal regular model of the elliptic curve over the ring of algebraic
integers~$ \OO_F $ of the number field~$ F $, the multiples of $ P $ hit at most
two irreducible components in each fibre.
We formulate our result for a general regular, geometrically connected and projective curve~$ C $ over
a number field~$ F $, with regular, flat and proper model over~$ \OO_F $,
but the proof works over various other bases too.

We recall from~\cite[\S1]{LJ15} the subgroup~$ K_2^T(C) $
in~$ K_2(F(C)) $, defined using the tame symbol on~$ C $,
with its subgroup~$\INTC $ of integral elements in~$ K_2^T(C) $,
defined using the tame symbol on a regular proper
model~$ \cC $ over~$ \OO_F $ of~$ C $,
with~$ \INTC $ independent of the model~\cite[Proposition~4.1]{LJ15}.

Here we shall assume~$ \cC $ is also flat over~$ \OO_F $,
which places us in the context of \cite[Chapter~9]{Liu}.
Because of the regularity of~$ \cC $, this only excludes connected components
contained in closed fibres, which do not influence the definition
of~$ \INTC $.

For an element of~$ K_2^T(C) $, the obstruction to its integrality at the fibre $ \cC_\PP $ of~$ \cC $ at~$ \PP $ 
is in a finitely generated Abelian group, and is trivial except for finitely many~$ \PP $.
So if for some~$ \a $ in~$ K_2^T(C) $ this obstruction is torsion at all fibres, then~$ M_\a \a $ is in~$ \INTC $
for some positive integer~$ M_\a $ that depends on~$ \a $
and over which we have no control.
But because the obstruction at a regular fibre~$ \cC_\PP $ is
always an element of~$ (\OO_F/\PP)^\times $,
by~\cite[Corollary~16.2]{mil71}
one can find~$ \b_\a $ in~$ K_2(F) $, unique up to adding an
element in the finite group~$ K_2(\OO_F) $, such that the obstruction
of $ \a + \b_\a $ is trivial for all regular fibres, and unchanged
for the finitely many singular fibres.
Its obstruction to integrality is now in a finitely generated Abelian
group, and because we assumed it was torsion, 
there is a positive integer~$ M $, depending only on~$ \cC $,
such that~$ M \a + M \b_\a $ is in~$ \INTC $. Its image under the regulator map
is~$ M $ times that of~$ \a $, as opposed to~$ M_\a $ times it
above. Thus, using~$ K_2(F) $ can help to find a larger subgroup
of~$ \INTC $.

For computing the tame symbol on $ \cC $ for any element of~$ K_2^T(C) $, we
may consider one fibre~$ \cF = \cC_\PP $ at~$ \PP $
at a time. For this we may and do replace~$ \OO_F $ with its localisation at~$ \PP $ and~$ \PP $
with the resulting maximal ideal. Let~$ \pi $ be a uniformising parameter
for~$ \PP $ and~$ k $ the residue field.

In order to compute the tame symbol $ T_\DD(\a) $ at every irreducible component~$ \DD $
of~$ \cF $, i.e., a curve in~$ \cF $ that is irreducible over~$ k $,
we may consider the divisor~$ (f) $ on~$ \cC $ for each of the functions~$ (f) $ involved in~$ \a $.
Then~$ (f) = (f)_h + (f)_v $ with~$ (f)_h $ the `horizontal part',
consisting of the divisor
of~$ f $ on $ C $ with the closed points of~$ C $ giving irreducible
curves on~$ \cC $  (with $ F $-rational points corresponding to sections
of~$ \cC $), and~$ (f)_v $ the `vertical part', a sum
over the irreducible components of~$ \cF $ with integer coefficients.
The properties of the intersection product on~$ \cC $
(see \cite[Corollary~8.3.6(b), Theorem~9.1.23]{Liu})
ensure~$ (f)_h $ determines $ (f)_v $ uniquely up to a multiple
of~$ [\cF] = (\pi) = \sum_\DD m_\DD [\DD] $, where the
sum runs through all the irreducible components of~$ \cF $ and
the~$ m_\DD $ are positive integers.
Because~$ f $ is determined by~$ (f)_h $ up to multiplication
by~$ F^\times $, and by~$ (f) $ up to multiplication
by~$ \OO_F^\times $,
such a computation can in general determine~$ T_\DD(\a) $ 
only up to multiplication by the image of~$ \OO_F^\times $ on~$ \DD $.
More information is needed for a more precise result.

The following (detailed) result, which seems the first general result not
ignoring any torsion, covers the situation needed in the proof of Theorem~\ref{thm:integrality}.
Its proof is based on the fact that
for at least one of the functions~$ f $, $ g $ and~$ h = - g/f $,
the points in its divisor hit only one irreducible component of~$ \cF $.
Because $ k^\times $ is finite, the obstruction
to integrality for the element~$ \a $ below is always of finite exponent, but
for an elliptic curve we get a more precise result in Corollary~\ref{new-int-cor}(2).

\begin{prop} \label{new-int-prop}
Let $ F $, $ \OO_F $, $ k $, $ C $, $ \cC $ and~$ \cF $ be as above.
Suppose we have non-zero functions $ f $ and~$ g $ on~$ C $,
distinct $ F $-rational points~$ O, P $ and $ Q $ on~$ C $,
and a positive integer~$ N $, such that~$ (f) = N (P) - N (O) $, $ (g) = N (Q) - N (O) $,
and~$ f(Q) = g(P) = 1 $.
Then~$ \a = \{ f , g \} $ is in $ K_2^T(C) $,
the function~$ h = - g/f $ is defined at $ O $, and~$ h(O) $ is an $ N $th root of unity~$ \v $ in~$ \OO_F^\times $.

Let $ \bar O $ be the point in the intersection of~$ \cF $ and~$ O $, and similarly for $ P $ and~$ Q $.
If~$ \bar O $, $ \bar P $ and $ \bar Q $ are on at most two irreducible components
of~$ \cF $, then~$ T_\DD(\a) $ is a constant function on~$ \DD $
for every irreducible component~$ \DD $ of~$ \cF $. In fact, with~$ M $
the order of the image $ \bar\v $ of~$ \v $ in $ k^\times $,
the following hold.

(1)
If $ \bar O $, $ \bar P $ and~$ \bar Q $ are on one irreducible component,
then~$ T_\DD( \a ) = 1 $ for all~$ \DD $.

(2)
Suppose that~$ \bar O $, $ \bar P $ and
$ \bar Q $ lie on two irreducible components~$ \cA \ne \cB $ of~$ \cF $,
with $ \bar O $ on~$ \cA $.
Let~$ \cA^o $ be the connected part of~$ \cup_{\DD \ne \cB} \DD $ that contains~$ \cA $,
and $ \cB^o $ the connected part of $ \cup_{\DD \ne \cA} \DD $ that contains~$\cB$.
\begin{itemize}
\item
If one of $ \bar P $ and $ \bar Q $ is on $ \cA $ but the other
one on $ \cB $, then~$ T_{\DD}(\a) = 1 $ for all~$ \DD $ equal to $ \cA $ or contained in~$ \cB^o $.
Moreover, if the two points on $ \cA $ coincide then $ T_\DD(\a) = 1 $ for all~$ \DD $ in $ \cF $.

\item
If $ \bar P $ and~$ \bar Q $ are on $ \cB $ then~$ M $ is an exponent 
for $ T_\DD( \a ) $ for all $ \DD$ equal to $ \cB $ or contained in~$ \cA^o $.
Moreover, if $ \bar P = \bar Q $ then this holds for all $ \DD $
in~$ \cF $.
\end{itemize}
\end{prop}

\begin{proof}
The product formula implies that~$ \a $ is in~$ K_2^T(C) $.
Then~$ (h) = N (Q) - N(P) $, so $ h $ is defined at~$ O $,
and from~$ \a = \{ f , g \} = \{ f , h \} $ we find~$ 1 = T_O(\a) = h(O)^N $.

(1) Note that~$ \cF $ is connected by \cite[Corollary~8.3.6(b)]{Liu}.
We treat various cases.

\mybullet
$ \bar O $, $ \bar P $ and $ \bar Q $ distinct.
Because $ f $ is regular and non-zero at $ \bar Q $, we see that~$ (f)_v = 0 $.
Similarly $ (g)_v = 0 $ as well, hence~$ T_\DD(\a) = 1 $ for all $ \DD $ in~$ \cF $.

\mybullet
$ \bar P \ne \bar Q $ but $ \bar O = \bar P $ or $ \bar O = \bar Q $.
By anti-symmetry of~$ \a $, we may assume~$ \bar O = \bar P \ne \bar Q $.
Then~$ (f)_v = 0 $, and $ f_{|\cF} = 1 $ because its poles and zeroes cancel at $ \bar  O = \bar P $,~$ f(\bar Q) = 1 $,
and $ \cF $ is connected. Therefore~$ T_\DD(\a) = 1 $ for every $ \DD $ in~$ \cF $.

\mybullet
$ \bar O $ distinct from~$ \bar P = \bar Q $.
We write~$ \a = \{ f, h \} $. Now $ (h)_v = 0 $, and as in the
first case we find from~$ \bar P = \bar Q $ and $ h(\bar O) = \bar\v $
that~$ h_{|\cF} = \bar\v $.
Note that $ f $ is uniquely determined by its divisor on~$ C $
and $ f(Q) = 1 $. We can construct it as follows.
Let $ \tilde f $ in~$ F(C)^\times $ be such that~$ (\tilde f)_h = N (P) - N (O) $.
Note that~$ m_\cA = 1 $ in~$ [\cF] $ because~$ \cA $ is hit by
a section, so multiplying $ \tilde f $ with a suitable power of~$ \pi $,
we may assume that~$ (\tilde f) $ does not contain~$ [\cA] $.
Then~$ (\tilde f) = N (P) - N (O) $.
Because the model is regular at~$ \bar P $, the local ring there
is a unique factorisation domain,
so in it we can write~$ \tilde f = u z^N $ with $ z $ a local equation of $ P $
and~$ u $ a unit. Then~$ f = u(Q)^{-1} z(Q)^{-N} \tilde f $,
with~$ u(Q) $ in~$ \OO_F^\times $. Hence~$ (f)_v = -N (z(Q)) $,
and~$ T_\DD( \a ) = 1 $ for every $ \DD $ in~$ \cF $ because $ \bar \v^N = 1 $.

\mybullet
$ \bar O = \bar P =  \bar Q $.
We blow up this point, obtaining a new regular model with one extra
irreducible component in the new fibre, a copy of $ \P_k^1 $, which
contains the new~$ \bar O $, $ \bar P $ and~$ \bar Q $. If all coincide, we blow
up this new point.
After finitely many steps (cf.~the proof of \cite[Theorem~9.2.26]{Liu})
we find a chain of new irreducible components,
such that in the last irreducible component of the chain the three points~$ \bar O $,
$ \bar P $ and $ \bar Q $ do not all coincide.
Then~$ T_\DD (\a) = 1 $ for all irreducible components~$ \DD $ of the new fibre
by the earlier cases, among which are those of the original fibre~$\cF$.

(2)
By anti-symmetry of $ \a $ for $ P $ and~$ Q $, we may assume $ \bar Q $ is on~$ \cB $.

\mybullet
$ \bar O $ and $ \bar P $ distinct on $\cA$.
From $ f(Q) = 1 $ we find that $ (f)_v = 0 $, so that~$ f $ restricts to a
constant function on every $ \DD \ne \cA $, and $ f_{|\cB^o} = 1 $.
Therefore $ T_\DD(\a) $ is a constant function for every $ \DD \ne \cA $
in~$ \cF $.
From $ g(P) = 1 $ and $ \bar P \ne \bar O $ we see that $ [\cA] $ in $ (g)_v $
has coefficient~0.
Hence~ $ T_\DD(\a) = 1 $ if~$ \DD = \cA $ or contained in~$ \cB^o $.

\mybullet
$ \bar O = \bar P $ on $\cA$.
From $ f(Q) = 1 $ we find~$ (f)_v = 0 $ and~$ f_{|\cF} = 1 $.
Hence $ T_\DD( \a ) = 1 $ for every $ \DD $ in~$ \cF $.

\mybullet
$ \bar P $ and $ \bar Q $ distinct on $ \cB$.
We write $ \a = \{ f , h \} $ with~$ (h)_v = 0 $.
From $ h(O) = \v $ we find that $ (h)_v = 0 $, so that~$ h $ restricts to a
constant function on every $ \DD \ne \cB $, and~$ h_{|\cA^o} = \v $.
Therefore $ T_\DD(\a) $ is a constant function for every~$ \DD \ne \cB $
in~$ \cF $, and~$ M $ an exponent for it if $ \DD $ is contained
in~$ \cA^0 $.
From $ f(Q) = 1 $ and $ \bar Q \ne \bar P $ we see that~$ [\cB] $ in $ (f)_v $
has coefficient~0, so that~$ T_\cB(\a) = 1 $.

\mybullet
$ \bar P = \bar Q $ on $ \cB$.
Here we have $ \a = \{ f , h \} $ with~$ (h)_v = 0 $, and~$ h_{|\cF} = \bar \v $.
Therefore $ M $ is an exponent for~$ T_\DD( \a ) $
for all $ \DD $ in $ \cF $.
\end{proof}

\begin{corollary} \label{new-int-cor}
Suppose that~$ C $ is a regular, geometrically connected and projective curve over a number field~$ F $,
with $ F $-rational points $ O $, $ P $ and $ Q$ and positive integer~$ N $ as in Proposition~\ref{new-int-prop}.
Suppose that for some regular, flat and proper model~$ \cC $ over~$ \OO_F $ of $ C $,
for each maximal ideal~$ \PP $ of~$ \OO_F $, 
$ O $, $ P $ and~$ Q $ hit at most two irreducible components in the fibre of~$ \cC $
over~$ \PP $.
Let~$ f $, $ g $, $ h = - g /f $ and~$ \a $ be as in the proposition.
Let~$ M' $ be the order of $ h(O) $ in $ \OO_F^\times $, 
which divides both $ N $ and the number of roots of unity in~$ F $.
Then the following hold.
\begin{enumerate}
\item
If for every~$ \PP $  the sections corresponding to~$ O $,
$ P $ and~$ Q $ hit the same irreducible component in the
fibre over~$ \PP $, then~$ \a $ is in~$ \INTC $.

\item
If the curve is an elliptic curve~$ E $, then $ M' \a $ is in~$ \INTE $.
\end{enumerate}
\end{corollary}

\begin{proof}
(1)
This is immediate from Proposition~\ref{new-int-prop}~(1).

(2)
For an elliptic curve~$ E $
there is a minimal flat, regular, proper model~$ \cE $, and every other
flat, regular proper model maps to it. On~$ \cE $ also at most two irreducible components can
be hit by the sections given by~$ O $, $ P $ and~$ Q $.
Moreover,
from the list of possible fibres~$ \cF $
\cite[p.486]{Liu},
one sees that $ \cF = \cA \cup \cB^o = \cA^o \cup \cB $ for all geometrically
irreducible components
$ \cA \ne \cB $ of multiplicity~1 (which are those that
can be hit by a section).
So for~$ \cE $, every irreducible component~$ \DD $
in any fibre over some~$ \PP $ is covered by Proposition~\ref{new-int-prop}, and~$ T_\DD( M' \a) = 1 $ always.
\end{proof}

\begin{remark}
Note that in part~(2) of the corollary, one only needs to multiply
by the least common multiple~$ M'' $ of the order of those~$ M = M_\PP $
for the~$ \PP $ where part~(2) of the proposition is required.
If all such~$ \PP $ have the same residue characteristic~$ p $, then
$ M'' $ is obtained from $ M' $ by dividing out all its factors~$ p $.
\end{remark}

The next example shows that in the situation of~ Proposition~\ref{new-int-prop}(2)
there does not always exist~$ \b $ in~$ K_2(F) $ with~$ T_\DD(\a) = T_\DD(\b) $ for all $ \DD $ in $ \cF $.
In the example~$ T_\DD(\a) $ depends on~$ \DD $, but~$ [\cF] = [\cA] + [\cB] $,
so~$ T_\cA(\b) = T_\cB(\b) $ for any~$ \b $ in~$ K_2(F) $.

\begin{example}
On the elliptic curve over~$ \Q $ defined by $ y^2 = x^3 + 1 $,
we let~$ P = (2,3) $ and~$ Q = - P = (2, -3) $.
Then $ (y+1) = 3 (-2P) - 3 (O) $ with~$ - 2 P = (0,-1) $, and~$ (y-2x+1) = 2 (P) + (-2P) - 3 (O) $, so we can take, with $ N = 6 $,
\begin{alignat*}{2}
f & = \frac1{108} \frac{(y-2x+1)^3}{y+1}
\qquad
g & = \frac1{108} \frac{(-y-2x+1)^3}{-y+1}
.
\end{alignat*}
Then $ h = -g/f = \frac{(-y-2x+1)^3}{(y-2x+1)^3} \frac{y+1}{y-1} $ satisfies $ h(O) = - 1 $.
The reduction of the curve at $ p =3 $ is of type~IV, so has only irreducible components
$ \cA$ and $ \cB$, with $ - [\cA] \cdot [\cA] = - [\cB] \cdot [\cB] =[\cA] \cdot [\cB] = 2 $. One has $ \bar P \ne \bar Q $
but both are on~$ \cB $.
On the flat minimal regular model over $ \Z_{(3)} $ one has $ (f) = 6 (P) - 6 (O) - 3 \cA $
and~$ (g) = 6 (Q) - 6 (O) - 3 \cA $, so~$ T_\cA (\{ f, g\} ) = -1 $
and~$ T_\cB (\{ f, g\} ) = 1 $.
\end{example}

\begin{remark} \label{remark:fafb}
One can employ the techniques used to prove Proposition~\ref{new-int-prop}
also in other ways. As an example, suppose that $ P $ is
an $ F $-rational point of prime order~$ p $ on an elliptic curve~$ E $
over a number field~$ F $, and that~$ P $ hits the 0-component
in every fibre of the minimal regular model~$ \cE $ of~$ E $ over~$ \OO_F $.
For a maximal ideal~$ \PP $ of~$ \OO_F $, the reduction~$ \bar P $ of~$ P $
in the (non-singular part of the) fibre~$ \cF_\PP $ of~$ \cE $ has order~1 or~$ p $.
If the order is~$ p $, then $ O, P, \dots, (p-1) P $ all have
different reductions at~$ \PP $. If the order is~1, then
translation by~$ P $ on the genus one curve~$ E $ induces
an automorphism of~$ \cE $ of order~$ p $ that maps~$ \bar O $
to~$ \bar P = \bar O $.
Hence it induces an automorphism of order~$ p $ of the blowup of~$ \cE $
at the closed point~$ \bar O $ in~$ \cF_\PP $. Because of the
action of this automorphism, the reductions of~$ O, P, \dots, (p-1) P $
in the new fibre above~$ \PP $ are either all different (but in the same
irreducible component) or all the same.
If they are the same, we repeat the procedure, obtaining from
the automorphism of the blowup an automorphism of order~$ p $ of the second blowup.
Repeating this
as necessary for all~$ \PP $ with~$ \bar P =  \bar O $
in~$ \cF_\PP $, we find a regular, flat and proper model~$ \cE' $ of~$ E $ such
that in every fibre~$ \cF_\PP' $ the reductions of~$ O, P, \dots, (p-1)P $
are distinct but in the same irreducible component.

For a subset~$ Z $ of closed points of~$ E $, let us say that
an element of~$ K_2(F(E)) $ is supported in~$ Z $ if 
it can be written using symbols with as entries functions with divisors supported in~$ Z $.
We want to describe the subgroup~$ \W $ of~$ K_2^T(E) $ supported in~$ O, P, \dots, (p-1) P $
as well as its intersection $ \Wint $ with $ \INTE $.
For this, we assume~$ p \ge 5 $, and introduce some notation.

For~$ a = 1,\dots,p-1 $, let~$ L^a $ be the subgroup of elements
in~$ K_2(F(E)) $ supported in~$O, P, 2P, \dots, a P $, so that $ \W $
consists of the elements of $ L^{p-1} $ with trivial symbol at
all multiples of~$ P $.
For~$ b\ne c
$ in~$ \{ 2,\dots,p-1\} $, let~$ f_{b,c} $ in~$ F(E)^\times $
be the function with divisor~$ (bP)- b (P) + (b-1) (O) $ on~$ E $, and~$ f_{b,c}(cP) = 1 $.
Because of the way the sections given by~$ O, P, \dots, (p-1) P $
lie on our model, its divisor on~$ \cE' $ is also~$ (bP)- b (P) + (b-1) (O) $.
For $ b = 2, \dots, p-1 $ we let~$ g_b $ in $ F(E)^\times $
be the function with divisor~$ p(P) - p(O) $ on~$ E $, and~$ g_b(bP)=1 $.
Then its divisor on~$ \cE' $ is also~$ p(P) - p(O) $.
The tame symbol of each~$ \{f_{b,c}, f_{c,b}\} $
or~$ \{f_{b,c} , g_b \} $ is trivial for every irreducible curve~$ \DD $ in $ \cE' $ except possibly
if~$ \DD $ is the section given by~$ P $ or~$ O $.
From the divisors of the functions on $ \cE' $, we see
that the values there are in~$ \OO_F^\times $.

For $ b = 2, \dots, p-1 $, a function in~$ F(E)^\times $
with divisor supported in~$ O, P, \dots, b P $
is a product of a power of~$ f_{b,d} $ for some fixed $ d \ne b $ and a
function with divisor supported in~$ O, P, \dots, (b-1) P $.
So if~$ \a $ is in $ L^b $, then modulo~$L^{b-1} $
it can be written as~$ \{ f_{b,d}, h \} $ with~$ (h) $
supported in~$ O, P, \dots, (b-1) P $. Therefore $ T_{bP}(\a) = 1 $
is equivalent to~$ h(bP) = 1 $, which means that~$ h $ is a product of
the~$ f_{c,b}^{\pm1} $ with $  2 \le c < b $ as well as $ g_b $.
Because~$ f_{b,d} $ scales to each~$ f_{b,c} $,
we see that~$ \a $ modulo~$ L^{b-1} $ is a sum of terms~$ \a_b $ of~$ \pm \{f_{b,c}, f_{c,b} \} $,
with $ 2 \le c < b $, and~$ \pm \{f_{b,c}, g_b \} $.

If we start with $ \a $ in~$ L^{p-1} $ and with trivial tame
symbol at $ 2 P, \dots, (p-1) P $, then~$ \a - \a_{p-1} $ is
in~$ L^{p-2} $, and again has trivial tame symbol at~$ 2 P, \dots, (p-1) P $.
Continuing this way, 
we can write $ \a = \a_{p-1} + \a_{p-2} + \dots + \a_2 + \a' $
with $ \a_b $ for $ b = 2,\dots,p-1 $ as above, and $ \a' $ 
in~$ L^1 $.
Then~$ \a'= \{ g_2, u \} + \g $ with~$ u $
in~$ F^\times $ and~$ \g $ in $ K_2(F) $,
and~$ \a $ is in~$ \W $ if and only 
if~$ T_P(\a_2+\dots+\a_{p-1}) = u^p $.
But~$ T_P(\a_j) $ for $ j = 2,\dots,p-1 $ is in~$ \OO_F^\times $,
so that~$ u $ must then be in~$ \OO_F^\times $. From the divisor
of~$ g_2 $ on~$ \cE' $, we see that then~$ \a_2+\dots+\a_{p-1} + \{g_2, u\} $
is in~$ \W \cap \INTE = \Wint $, so that $ \W = \Wint + K_2(F) $,
and $ \a $ is in $ \Wint $ if and only if~$ \g $ is in~$ K_2(\OO_F) $.
It also follows that the obstruction to completing any sum~$ \b $ of
terms~$ \pm \{f_{b,c}, f_{c,b} \} $
and~$ \pm \{f_{b,c}, g_b \} $ to an element of $ \INTE $ by adding
an element of $ L^1 $ to it lies in 
in~$ \OO_F^\times / (\OO_F^\times)^p $.
If this is possible, then the $ \b + \{g_2, u\} + \g $ with
$ u^p = T_P(\b) $ in $ \OO_F^\times $
and~$ \g $ in~$ K_2(\OO_F) $ are all such completions to an element
of~$ \INTE $, and this way we obtain all elements of~$ \Wint $.
In particular,  with~$ u_{b,c} = T_P(\{ f_{b,c}, f_{c,b} \}) $
in~$ \OO_F^\times $, the element~$ \b_{P,b,c} = p \{ f_{b,c}, f_{c,b} \} + \{g_2 , u_{b,c} \} $
is in~$ \Wint \subseteq \INTE $.
\end{remark}

\section{\texorpdfstring{$K_2$}{K2} of families of elliptic curves with 7, 8 and 10-torsion points} \label{K2-families}

In this section, we give our main construction of families of elliptic curves $E$ over (in
general non-Abelian) cubic and quartic fields, with respectively three and four elements in~$\INTE$
that are, in general, linearly independent, as stated in Theorem~\ref{thm:independence}.
(See Section~\ref{K2-cubic} for a simpler construction and proof
for cubic fields.)

Our starting point is Bloch's construction of elements in~$K_2^T(E)$
using the torsion subgroup of $E$, which we recall and refine
in Section~\ref{blochsection}.
To obtain elements in~$ K_2^T $ of curves in this way, it is natural
to consider elliptic curves with a large enough rational torsion
subgroup. For example, Liu and de Jeu~\cite[Proposition 6.9]{LJ15}
construct families of elliptic curves over quadratic fields containing
a subgroup isomorphic to~$\Z/2\Z \times \Z/4\Z$,
with two elements in~$ \INTE $ that are in general linearly independent.

In Section~\ref{subsection:parametrization} we apply
Bloch's construction to the universal elliptic
curve $\mathcal{E}_1(N)$ over the modular curve $Y_1(N)$. This
is the universal family of elliptic curves with a point~$ P $ of order $N$.
The definitions of $Y_1(N)$ and $\mathcal{E}_1(N)$ are reviewed in
Section~\ref{subsection:parametrization}. We focus here on the cases
$N=7,8,10$, for which the modular curve $X_1(N)$ has genus~$0$.
In particular, we work out the equations of these families over $\Q(t)$,
where~$t$ is a Hauptmodul of $X_1(N)$ expressed in terms of the
Weierstra\ss{} $\wp$-function.

In fact, one can show for any $N$ that the element~$ S_{P,s} $
we construct in $K_2^T$ of a given fibre $\mathcal{E}_1(N)_t$ of $\mathcal{E}_1(N)$
is, after tensoring with~$ \Q $, a non-zero rational multiple of
Beilinson's Eisenstein symbol associated to the point $sP$ \cite{Be86, DS91, DK11}.

In Section~\ref{subsec:integrality}, we find sufficient conditions on
the parameter $t$, belonging to a given number field, so that
the elements from Bloch's construction are in $\INTE$. This uses
the new integrality criterion from Section~\ref{new-integrality}.
The conditions found lead naturally to the construction of our
cubic and quartic fields in Lemma~\ref{lemma:unit}.

The main theorem of this section is then formulated in Section~\ref{subsec main result}.
That the elements we constructed are in~$ \INTE $ has by then already
been established, but not yet their linear independence.
This independence follows from a stronger result on the limit behaviour of the
Beilinson regulator (see the definition in Section~\ref{conjecture}). The proof
of this is somewhat involved, and is given in Section~\ref{section:independence}.

\subsection{Bloch's construction of elements in \texorpdfstring{$\textit K_{\mathbf 2}^{\mathbf T}$}{K2T}} \label{blochsection}
We revisit and refine this construction from~\cite[Section 10.1]{Bl00}.
Let $E$ be an elliptic curve defined over a field~$F$, and~$P$ an~$F$-rational point
on~$ E $ of order $N$. For $1 \leqslant s \leqslant N-1$, let $ f_{P,s}$
in~$F(E)^\times$ be a function with divisor $(f_{P,s}) = N(sP) - N(O)$.
Then for $ s \ne t $, the symbol $ T_{P,s,t} = \{ f_{P,s} / f_{P,s}(tP) , f_{P,t} / f_{P,t}(sP) \} $
is in $ K_2^T(E) $ and independent of the choice of the functions
(cf.~\cite[Construction~4.1]{DJZ}).
For $ 1 \le s \le N-1 $ we let
$S_{P,s} = \sum_{t=1, t \ne s}^{N-1} T_{P, s,t} $.
Adding the trivial symbol~$ \{f_{P,s}, - f_{P,s} \} $
we see that, modulo symbols with (at least) one function in~$ F^\times$,
we have~$ S_{P,s} \equiv \{ f_{P,s}, g_P \} $
for some~$ g_P $ in~$ F(E)^\times $ with~$(g_P) = \sum_{t=1}^{N-1} N(tP) - N(N-1) (O) $.

For an integer~$ a $ with~$ \gcd(a,N) = 1$, $ a P $ also has
order~$ N $, and~$ T_{aP, s,t} = T_{P, as, at} $ because
the elements are uniquely determined by the divisors of the functions
involved. We may consider $ s $ and $ t $ as elements of~$ \Z/N\Z $,
which gives
\begin{equation*}
S_{aP,s} = \sum_{t \ne 0, s \in \Z/N\Z} T_{aP, s, t} = \sum_{t \ne 0, s \in \Z/N\Z} T_{P, as, at} = S_{P,as}
\end{equation*}
as multiplication by~$ a $ permutes $ \Z/N\Z $.
So the set $ \{ S_{P,1}, \ldots, S_{P,N-1} \} $
depends only the subgroup generated by~$ P $ and not on the choice
of a generator.

We recall that the action of the correspondences on $ E \times E $ on $ K_2^T(E)/K_2(F) $
implies that pulling back along multiplication by~$ -1 $
on~$ E $ induces multiplication by $ -1 $ on $ K_2^T(E) / K_2(F) $.
Therefore, the~$ S_{P,s} $ with~$ 1 \le 2s \le N-1 $ will generate
the same subgroup of~$ K_2^T(E) $ modulo~$ K_2(F) $ as all~$ S_{P,s} $.

\subsection{Families of elliptic curves with torsion} \label{subsection:parametrization}

We first review the definition of the modular curve $Y_1(N)$
over $\Q$. For $N \geqslant 4$, we consider the functor
sending a $\Q$-scheme $S$ to the set of isomorphism
classes of pairs $(E, P)$ where $E$ is an elliptic curve
over $S$, and $P \in E(S)$ is a point of exact order $N$.
This functor is representable by a $\Q$-scheme which
we denote by $Y_1(N)$ (see \cite[Theorem 8.2.2]{DI95}).
Moreover, there is a universal elliptic curve $\mathcal{E}_1(N)$ over
$Y_1(N)$, equipped with a section of order $N$.

We shall need a description of the complex points of $Y_1(N)$
and $\mathcal{E}_1(N)$.
Let $\H$ denote the
upper half-plane. By \cite[Sections 1 and 2]{Kat04},
there is an isomorphism of Riemann surfaces
\begin{equation*}
\begin{split}
\nu\colon\Gamma_1(N) \backslash \H & \to Y_1(N)(\C) \\
\tau & \mapsto \Bigl(\C/(\Z+\tau\Z), \Bigl[\frac{1}{N}\Bigr]\Bigr),
\end{split}
\end{equation*}
where
\begin{equation*}
\Gamma_1(N) = \Bigl\{A \in \mathrm{SL}_2(\Z) : A \equiv
\begin{pmatrix} 1 & * \\ 0 & 1 \end{pmatrix} \bmod{N}\Bigr\}
\end{equation*}
acts on $\H$ by $(\begin{smallmatrix} a & b \\ c & d \end{smallmatrix})
\cdot \tau = \frac{a\tau+b}{c\tau+d}$. Moreover, the fibre of $\mathcal{E}_1(N)(\C)$
over the point $\nu(\tau)$ can be identified with $\C/(\Z+\tau\Z)$. In fact, $\mathcal{E}_1(N)(\C)$ can be described as a quotient of $\H \times \C$ by the semi-direct product $\Gamma_1(N) \rtimes \Z^2$ \cite[7.1.1]{DK11}. For any $z \in \C$, we denote by $[\tau, z]$ the point
of $\mathcal{E}_1(N)(\C)$ corresponding to $z \in \C/(\Z+\tau\Z)$. From the moduli
description, we deduce that for a matrix
$A = (\begin{smallmatrix} s & t \\ u & v \end{smallmatrix})$ in $\Gamma_1(N)$,
we have $[A\tau, z] = [\tau, (u\tau+v) z]$.

In this article, we shall use an alternative description of
$Y_1(N)(\C)$. Let $ W_N $ be the Atkin-Lehner involution
on $\Gamma_1(N) \backslash \H$, induced by $W_N(\t) = - \frac{1}{N\t}$.
Consider the composite isomorphism
\begin{equation*}
\nu' = \nu \circ W_N\colon\Gamma_1(N) \backslash \H \to Y_1(N)(\C).
\end{equation*}
For $\t \in \H$, we have
\begin{equation*}
\nu'(\t) = \nu(W_N(\tau)) = \Bigl(\frac{\C}{\Z+W_N(\tau)\Z}, \Bigl[\frac{1}{N}\Bigr]\Bigr) \cong \Bigl(\frac{\C}{\Z+N\tau\Z}, [\tau]\Bigr),
\end{equation*}
where the last isomorphism is induced by multiplication by $N\tau$.
In this way, the fibre of $\mathcal{E}_1(N)(\C)$ over $\nu'(\tau)$ can be identified
with the elliptic curve $\C/(\Z+N\tau\Z)$.

\begin{remark} \label{remark nu}
This unusual description of $Y_1(N)(\C)$ with the map $\nu'$ is related
to considering another model of the modular curve, obtained by formulating
the moduli problem using embeddings $\mu_N \to E$,
where $\mu_N$ is the $\Q$-scheme of $N$-th roots of unity \cite[8.2.2]{DI95}.
In this article, we stick with the usual model of $Y_1(N)$ over~$\Q$ defined
at the beginning of the section, but work with $\nu'$ to describe its complex
points. One advantage of this description is that a modular function
$t \in \C(Y_1(N))$ is defined over a subfield $K \subset \C$ if and only if
the Fourier expansion of $t \circ \nu'$ at the cusp $i\infty$
has coefficients in $K$ \cite[Remark 2.15]{Str15}.
\end{remark}

Under the isomorphism $\nu\colon\Gamma_1(N) \backslash \H \to Y_1(N)(\C)$,
the complex conjugation on $Y_1(N)(\C)$ corresponds to the standard
complex conjugation on $\Gamma_1(N) \backslash \H$ induced by
$c(\tau) = -\bar{\tau}$. This follows from the algebraic description of the
elliptic curve $\C/(\Z+\tau\Z)$ using the Weierstra\ss\ parametrisation, and
the identity $\overline{\wp_\tau(z)} = \wp_{-\bar\tau}(\bar z)$ for the
Weierstra\ss{} $\wp$-function. For the same reason, the complex conjugation
on~$\mathcal{E}_1(N)(\C)$ sends $[\tau, z]$ to $[-\bar\tau, \bar z]$.
Moreover, since $W_N \circ c = c \circ W_N$ on $\H$, the isomorphism
$\nu'$ is also compatible with complex conjugation.

We now discuss families of elliptic curves with torsion, and their
complex-analytic parametrisations.
Let $F$ be a field and $N \geqslant 4$ an integer.
Every pair $(E, P)$ consisting of an elliptic curve $E$ over $F$
and an $F$-rational point $P$ of order $N$ admits a unique Weierstra\ss{} model
\begin{equation}\label{eqn:TNF}
E: y^2+(1-g)xy-fy = x^3-fx^2
\end{equation}
with $f$ in $ F^{\times}$, $g$ in $ F$, and where $P = (0,0)$.
This is called the \emph{Tate normal form} of~$E$ \cite[Lemma 2.6]{Ba}.
By our construction in Section~\ref{blochsection}, we
have elements~$ S_{P,s} $ in~$ K_2^T(E) $ for~$ s = 1,\dots, \lfloor \frac{N-1}{2}\rfloor$.

In particular, consider the universal elliptic curve $\mathcal{E}_1(N)$ over the
function field $F = \Q(X_1(N))$ of the smooth compactification $X_1(N)$
of $Y_1(N)$. It has a Tate normal form \eqref{eqn:TNF}, where
$P=(0,0)$ is the universal section of order $N$.
For certain values of $N$, one can give
an explicit parametrisation of $\mathcal{E}_1(N)(\C)$ in terms of elliptic functions
using a method of Lecacheux (see \cite[Section 3]{BN} in the case $N=8$).
Recall that the fibre of $\mathcal{E}_1(N)(\C)$ over $\nu'(\tau)$
with $\tau \in \H$, is canonically isomorphic to $\C/(\Z+N\tau\Z)$.
Then the Tate normal form of $\mathcal{E}_1(N)(\C)$ has a parametrisation
\begin{equation} \label{param E1N}
\begin{split}
x & = u(\tau)^2 \wp_{N\tau}(z) + r(\tau) \\
y & = u(\tau)^3 \Bigl(\frac{\wp_{N\tau}'(z)}{2}\Bigr) + u(\tau)^2 s(\tau) \wp_{N\tau}(z) + v(\tau)
\end{split}
\end{equation}
with $u(\tau) \in \C^\times$, $r(\tau), s(\tau), v(\tau) \in \C$, and where
$\wp_{N\tau}$ denotes the Weierstra\ss{} $\wp$-function on $\C/(\Z+N\tau\Z)$.

From now on, we assume that~$N = 7, 8$ or~$ 10$, in which cases the modular curve $X_1(N)$
has genus~$0$. For these values of $N$, the Tate normal form of $\mathcal{E}_1(N)$ and its
discriminant $\Delta$ are given in Table~\ref{tab: tnf},
where $t$ stands for a generator of~$F = \Q(X_1(N))$
(see \cite[Section 2.3]{Ba} for how to derive these equations).

\begin{table}
\begin{tabular}{cccc}
  \hline\vph
  $N$ & $f$ & $g$ & $\Delta$  \\
  \hline\vph
  $7$ & $t^3-t^2$ & $t^2-t$ & $t^7(t-1)^7(t^3 - 8t^2 + 5t + 1)$\\
  \hline\vph
  $8$ & $2t^2-3t+1$ & $\frac{2t^2-3t+1}{t}$ & $t^{-4}(t-1)^8(2t-1)^4(8t^2 - 8t + 1)$\\  \hline\vph
  $10$ & $\frac{2t^5-3t^4+t^3}{(t^2-3t+1)^2}$ & $\frac{-2t^3+3t^2-t}{t^2-3t+1}$ & $t^{10}(t-1)^{10}(2t-1)^5(t^2-3t+1)^{-10}(4t^2-2t-1)$\\
  \hline
\end{tabular}
\vskip 7pt
\caption{\label{tab: tnf} Tate normal form in~\eqref{eqn:TNF} for $N=7, 8, 10$.}
\end{table}

In the case $N=7$, we have
\begin{equation*}
2P = (t^2(t-1), t^3(t-1)^2), \quad 3P = (t(t-1), t(t-1)^2),
\end{equation*}
which, using \eqref{param E1N}, gives the parametrisation
\begin{equation} \label{eq t 7}
t = \frac{x(2P)-x(P)}{x(3P)-x(P)} \qquad
t \circ \nu'(\tau) = \frac{\wp_{7\tau}(2\tau)-\wp_{7\tau}(\tau)}{\wp_{7\tau}(3\tau)-\wp_{7\tau}(\tau)}.
\end{equation}
Similarly, in the case $N=8$, we compute the coordinates of $2P$
and $3P$ in terms of $t$. This results in the parametrisation
\begin{equation*}
t \circ \nu'(\tau) = \frac{\wp_{8\tau}(2\tau)-\wp_{8\tau}(\tau)}{\wp_{8\tau}(3\tau)-\wp_{8\tau}(\tau)},
\end{equation*}
Finally, in the case $N=10$, the computation of $2P$ and $4P$ gives
\begin{equation*}
t \circ \nu'(\tau) = \frac{\wp_{10\tau}(2\tau)-\wp_{10\tau}(\tau)}{\wp_{10\tau}(4\tau)-\wp_{10\tau}(\tau)}
\end{equation*}

We note that for $N=7,8,10$, the function $t \circ \nu'$ is defined
and non-vanishing on~$\H$, since $\wp_\tau(z) = \wp_\tau(z')$ only for
$z = \pm z'$. Moreover, \cite[Theorem 6.2(a)]{Sil94} shows that the Fourier
expansion of $t \circ \nu'$ has coefficients in $\Q$. By Remark \ref{remark nu},
this means that $t$ is a modular unit on $Y_1(N)$, that is, belongs to
$\mathcal{O}(Y_1(N))^\times$.

We shall need the values of $t \circ \nu'$ at the cusps
for the proof of the main theorem in Section \ref{end proof}.
Especially relevant here are the cusps lying on the real
locus of $\Gamma_1(N) \backslash (\H \cup \Q \cup \{i\infty\})$,
by which we mean the fixed points under the complex
conjugation $c(\tau) = -\bar\tau$ defined above.
This real locus is completely determined by Snowden in
\cite[Section 6.4]{Sn}, for any $N$.

Let us first deal with the cusp $i\infty$.
Write, as usual,
$q = e^{2\pi i\tau}$. Using the Fourier expansion of the Weierstrass
$\wp$ function \cite[Theorem 6.2(a)]{Sil94}, we obtain
\begin{equation} \label{t at oo}
t \circ \nu'(\tau) = 1-q + O(q^2) \qquad (N = 7, 8, 10).
\end{equation}
In particular, $t \circ \nu'$ takes the value $1$ at the cusp $i\infty$.
The values at the other cusps are given by the poles of $j(E) \in \Q(t)$,
since the cusps of $X_1(N)$ lie over the cusp $i\infty$ of $X_1(1)$.
Moreover, the map $t \circ \nu'$ is compatible with complex conjugation,
thus induces a homeomorphism between the real locus of
$\Gamma_1(N) \backslash (\H \cup \Q \cup \{i\infty\})$ and $\P^1(\R)$.
Actually, all the cusps are real for these values of $N$.

Take for example $N=7$. We find that the images of the cusps
under $t \circ \nu'$ are $0, 1, \infty$ and the three real roots of
$t^3 - 8t^2 + 5t + 1$. By \cite[Section 6.4]{Sn}, the real locus of
$\Gamma_1(7) \backslash (\H \cup \Q \cup \{i\infty\})$ is the hyperbolic hexagon
whose vertices are the cusps $i\infty$, $\frac12$, $\frac37$,
$\frac13$, $\frac27$, $0$ in that order. We already know that
$t \circ \nu'(i\infty) = 1$, and more precisely
$t \circ \nu'(iy) \to 1^-$ when $y \to +\infty$. Since $t \circ \nu'$
preserves the cyclic order, this is sufficient to get the images
of all these cusps. The correspondence is given in Table
\ref{table cusps N=7}.

Similarly, for $N=8$, the real locus of
$\Gamma_1(8) \backslash (\H \cup \Q \cup \{i\infty\})$ is the
hyperbolic hexagon whose vertices are the cusps $i\infty$, $\frac12$,
$\frac38$, $\frac13$, $\frac14$, $0$ in that order. The images of these
cusps under $t \circ \nu'$ are given in Table \ref{table cusps N=8}.
Finally, for $N=10$, the real locus of
$\Gamma_1(10) \backslash (\H \cup \Q \cup \{i\infty\})$ is the
hyperbolic octagon whose vertices are the cusps $i\infty$, $\frac12$,
$\frac25$, $\frac13$, $\frac{3}{10}$, $\frac14$, $\frac15$, $0$ in that
order. The images of these cusps under $t \circ \nu'$ are given in
Table \ref{table cusps N=10}.

\begin{table}[t]
\begin{tabular}{c|cccccc}
  \hline\vph
  Cusps & $i\infty$ & $\frac12$ & $\frac37$ & $\frac13$&  $\frac27$ & $0$ \\
  \hline\vph
  $t \circ \nu'$ & $1$ & $7.295\dots$ & $\infty$ & $-0.158\dots$ & $0$ & $0.862\dots$ \\
  \hline
\end{tabular}
\vskip 7pt
\caption{\label{table cusps N=7} The cusps for $N=7$}
\end{table}

\begin{table}[t]
\begin{tabular}{c|cccccc}
  \hline\vph
  Cusps & $i\infty$ & $\frac12$ & $\frac38$ & $\frac13$ & $\frac14$ & $0$ \\
  \hline\vph
  $t \circ \nu'$ & $1$ & $\infty$ & $0$ & $\frac{2-\sqrt{2}}{4}$ & $\frac12$ & $\frac{2+\sqrt{2}}{4}$ \\
  \hline
\end{tabular}
\vskip 7pt
\caption{\label{table cusps N=8} The cusps for $N=8$}
\end{table}

\begin{table}[t]
\begin{tabular}{c|cccccccc}
  \hline\vph
  Cusps & $i\infty$ & $\frac12$ & $\frac25$ & $\frac13$ & $\frac3{10}$ & $\frac14$ & $\frac15$ & $0$ \\
  \hline\vph
  $t \circ \nu'$ & $1$ & $\frac{3+\sqrt{5}}{2}$ & $\infty$ & $\frac{1-\sqrt{5}}{4}$ & $0$ & $\frac{3-\sqrt{5}}{2}$ & $\frac12$ & $\frac{1+\sqrt{5}}{4}$ \\
  \hline
\end{tabular}
\vskip 7pt
\caption{\label{table cusps N=10} The cusps for $N=10$}
\end{table}

\begin{remark}\label{remark:isoN}
For $N=7$, let $t' = 1/(1-t)$. We have an isomorphism from $E_t$ to~$E_{t'}$
given by mapping~$ (x,y) $
to~$ \bigl( (t-1)^{-4} (x - t^2+ t ), (t-1)^{-6} (y - (t^2 - 2t) x - t^3 + t^2 ) \bigr) $.
With~$P'=(0,0)$ on $E_{t'}$, it is easy to check this maps $P=(0,0)$ on~$ E_t $ to $2P'$.
Similarly, $E_t$ is also isomorphic to $E_{t'}$ for~$t'=1-\frac{1}{t}$,
with~$P$ mapping to~$3P'$.

Also, for $N=8$, $ t $ and~$t' = 1-t$ give isomorphic elliptic curves,
and the same holds for $N=10$ for~$ t $ and~$t' = \frac{t-1}{2t-1}$,
with~$P$ mapping to $3P'$ in both cases.

Because under these isomorphisms, $ P $ is mapped to $ a P' $
for an integer~$ a $ with $ \gcd(a,N) = 1 $, so~$ S_{P,s} $
corresponds to~$ S_{aP',s} $. By what we saw in Section~\ref{blochsection},
the~$ S_{P,s} $ for $ s = 1,\dots,N-1 $ on~$ E_t $
and the~$ S_{P',s} $ for $ s = 1,\dots,N-1 $ on~$ E_{t'} $ correspond,
in some order.
\end{remark}

\subsection{Integrality of the elements} \label{subsec:integrality}
We can now, under suitable assumptions, relate the elements~$ S_{P,s} $,
as mentioned after~\eqref{eqn:TNF} for~$ N = 7 $, 8 and~10,
to~$ \INTE $.
It will be easier to use those in Section~\ref{section:independence} instead of the $ T_{P,s,t} $ (but see Remark~\ref{remark:ST}).

\begin{theorem}\label{thm:integrality}
Let $N =7, 8 $ or~$ 10$, and let $F$ be a number field.
For a fixed~$t$ in~$ F$ such that $\Delta(t) = \Delta_N(t)$ in~Table~\ref{tab: tnf}
is defined and non-zero, let $E$ be the elliptic curve over~$ F $
given by the corresponding Tate normal form~\eqref{eqn:TNF},
and $ P = (0,0) $.
Then, under the following assumptions, $ 2 P $ hits the 0-component
in each fibre of the minimal regular model of~$ E $ over~$ \OO_F $.
\begin{enumerate}
  \item $t, 1-t $ are in $\OO_F^{\times}$, for $ N = 7 $.
  \item $\frac{1}{t}-1, \frac{1}{t}-2 $ are in $\OO_F^{\times}$, for~$ N = 8 $.
  \item $\frac{1}{t}-1, 1-2t $ are in $ \OO_F^{\times}$, for~$ N = 10 $.
\end{enumerate}
Moreover, under those assumptions, the following hold for~$ s = 1,\dots, N-1 $.
\begin{itemize}
\item
For $ N = 7 $, $ S_{P,s} $ is in~$ \INTE $.

\item
For $ N = 8 $ and 10, $ N' S_{P,s} $ is in~$ \INTE $,
where~$ N' $ is the greatest common divisor of~$ N $ and the number
of roots of unity in~$ F $.
\end{itemize}
\end{theorem}

\begin{proof}
Below, we shall show that~$ 2 P $ always hits the 0-component,
so that every multiple of~$ P $ hits either that component or
the component hit by~$ P $.
Note that it follows that only the 0-component is hit if~$ N $ is odd.
We constructed~$ S_{P,s} $ in Section~\ref{blochsection}
as a sum of elements~$ T_{s,t} $, to which
part~(1) of Corollary~\ref{new-int-cor} applies if~$ N $ is odd, and
part~(2) if~$ N $ is even, and the same holds for their sum.
Therefore our statement about~$ \INTE $ follows from our claim
on~$ 2 P $.

We shall prove it
by showing that, for every~$ \PP $, for a minimal Weierstra\ss\ equation of~$ E $ at~$ \PP $
the point $P$ or $2P$ reduces to a regular point for
it modulo~$ \PP $.
We do this separately for the cases~$ N = 7 $, 8 or~10.

(1) It is straightforward to calculate the invariants $c_4$ and $\Delta$ of $E$ (see \cite{Ba})
\begin{align*}
c_4 &= (t^2-t+1)(t^6-11t^5+30t^4-15t^3-10t^2+5t+1),\\
%c_6 &= -t^{12} + 18t^{11} - 117t^{10} + 354t^9 - 570t^8 + 486t^7 - 273t^6 + 222t^5 - 174t^4 + 46t^3 + 15t^2 - 6t - 1,\\
\Delta &= t^7(t-1)^7(t^3-8t^2+5t+1).
\end{align*}
Since $t$ and $ 1-t$ are in $ \OO_F^{\times}$, we only consider places~$ \PP $ with~$\PP \mid (t^3-8t^2+5t+1)$.

By Lemma 3.1 of \cite{Ba} with $ r = 1 $ and $ s = 0 $, there are $\mu,\nu$ in $\OO_F$ such that $\mu c_4+\nu \Delta = 7^2$.
If $ \PP \nmid 7 $ then~$\PP \nmid c_4$, and the Weierstra\ss\ equation~\eqref{eqn:TNF} is minimal at~$\PP$.
The derivative with respect to $y$ at the reduction of~$P=(0,0)$ is
the reduction of~$f$, which is non-zero.

If $\PP \mid 7$, we transform the equation into
\begin{equation*}
y^2 = x^3 - 27c_4x - 54c_6
\end{equation*}
where $P$ is transformed into $(3(1-g)^2-12f, -108f)$
(see \cite[Ch III, Section 1]{Si} for the transformation).
If the equation is minimal, this point reduces to a regular point
since $\PP \nmid 108$ and $ f $ is in $ \OO_F^\times $.
If the equation is not minimal, we can transform
it into a minimal Weierstra\ss\ equation
$$y^2 = x^3 - 27\frac{c_4}{u^4}x - 54\frac{c_6}{u^6}$$
for certain~$ u $ with $v_{\PP }(u)>0$ and $P$ is transformed into $\left(\frac{3(1-g)^2-12f}{u^2}, -\frac{108f}{u^3}\right)$ which is reduced to the point at infinity.

(2) By assumption, $t' = \frac{1}{t}$ is in $ \OO_F$.
Under the change of variables $x\mapsto t^2 x, y\mapsto t^3 y$
the elliptic curve $E$ is isomorphic to
\begin{equation*}
E': y^2+(-t'^2+4t'-2)xy-t'(t'-1)(t'-2)y = x^3-(t'-1)(t'-2)x^2
\end{equation*}
with invariants
\begin{align*}
c_4 &= t'^8-16t'^7+96t'^6-288t'^5+480t'^4-448t'^3+224t'^2-64t'+16,\\
\Delta &= (t'-2)^4(t'-1)^8 t'^2 (t'^2-8t'+8) ,
\end{align*}
with $ t'-2 $ and $ t'-1 $ in $ \OO_F^\times $ by assumption.
So at a place~$ \PP $ with $v_\PP(t') > 0$, we have $v_\PP(2) = 0$, hence
$ v_\PP(c_4) = 0 $.
Therefore
the given Weierstra\ss\ equation of $E'$ is minimal at~$ \PP $.
The Weierstra\ss\ equation has reduction
$ y^2 -2 xy = x^3 -2 x^2$,
and the point $2P = ((t'-1)(t'-2),(t'-1)^2(t'-2)^2)$ reduces
to the regular point~$ (2,4) $.

At a place~$ \PP $ with $v_\PP(t'^2-8t'+8) > 0 $, we also have $v_\PP(2) = 0$
because $ v_\PP(t'-2) = 0 $.
By Lemma~3.1 of~\cite{Ba}, there exist $\mu,\nu$ in $\OO_F$ such that $\mu c_4 + \nu \Delta = 2^4$.
It follows that~$v_\PP(c_4)=0$, so that the above Weierstra\ss\
equation is minimal at~$ \PP $.
In this case $v_\PP(t') = 0$ as well,
and $P = (0,0)$ reduces to a regular point on the curve defined
by the reduction of the Weierstra\ss\ equation.

(3) By assumption, $t' = \frac{1}{t}$ is in $ \OO_F$, and $ t' -1 $
as well as $ 1 - 2 t'^{-1} = \frac{t'-2}{t'} $ are in~$ \OO_F^\times $.
Under the change of variables $x\mapsto \frac{x}{t'^2(t'^2-3t'+1)^2}, y\mapsto \frac{y}{t'^3(t'^2-3t'+1)^3}$
the elliptic curve $E$ is isomorphic to~$ E' $, defined by
\begin{equation*}
y^2+(t'^3-2t'^2-2t'+2)xy+t'^2(t'-1)(2-t')(t'^2-3t'+1)y = x^3+t'(t'-1)(2-t')x^2
,
\end{equation*}
with invariants
\begin{align*}
\begin{split}
  c_4 &= t'^{12} - 8t'^{11} + 16t'^{10} + 40t'^9 - 240t'^8 + 432t'^7 - 256t'^6 - 288t'^5 + 720t'^4 \\
  &\qquad\qquad\qquad\qquad\qquad\qquad\qquad\qquad\qquad\qquad - 720t'^3 + 416t'^2 - 128t' + 16,
\end{split}
\\
\Delta
  &= \left(\frac{t'-2}{t'}\right)^5(t'-1)^{10}t'^{10}(t'^2+2t'-4)(t'^2-3t'+1)^2.
\end{align*}
Suppose $\PP | \Delta$. We consider the two cases $\PP | 2$ and $\PP \nmid 2$.

\smallskip

\underline{Case 1: $\PP | 2$.}

First assume $\PP | t'$. We have $v_{\PP}(t')=v_{\PP}(t'-2)$ because
$\frac{t'-2}{t'} $ is in $ \OO_F^{\times}$, so that $v_{\PP}(t') \leqslant v_{\PP}(2)$.
Consider the change of variables $x \mapsto t'^2x$ and $y \mapsto t'^3y$ from $E'$ to
\begin{equation*}
E'': y^2+\frac{t'^3-2t'^2-2t'+2}{t'}xy+(t'-1)\frac{2-t'}{t'}(t'^2-3t'+1)y = x^3+(t'-1)\frac{2-t'}{t'}x^2
,
\end{equation*}
where the coefficients are in~$ \OO_F $.
The reduction of this equation modulo~$ \PP $ is
$ y^2 + (v+1) xy+v y = x^3+ v x^2 $
with $ v\ne 0$ the reduction of $ \frac{t'-2}{t'} $.

If $v_{\PP}(t') = v_{\PP}(2)$, then $v(c_4(E'')) =  v_{\PP}(c_4 / t'^{4}) = 0$.
Thus the Weierstra\ss\ equation of $E''$ above is minimal at $ \PP $.
The point $P = (0,0)$ does not reduce to a singular point.
On the other hand, if $v_{\PP}(t') < v_{\PP}(2)$, then $v_{\PP}(\Delta(E'')) = v_{\PP}(\Delta(E')/t'^{12}) = 0$
because $v_{\PP}(\Delta(E')=v_{\PP}(t'^{10})+v_{\PP}(t'^2+2t'-4) = v_{\PP}(t'^{12})$.
Thus $E''$ is minimal and has good reduction at~$\PP$.

If $\PP | (t'^2+2t'-4)$, then we have $\PP | t'$, which is already dealt with.

If $\PP | (t'^2-3t'+1)$, then $ v_\PP(t') = 0 $, hence $v_{\PP}(c_4) = 0$
also,
and the equation of~$E'$ is minimal at $\PP$.
The reduction
of~$2P = (t'(t'-1)(t'-2), -t'(t'-1)^2(t'-2)^2)$ on the curve
defined by the reduction of the Weierstra\ss\ equation is regular because the partial derivative
with respect to $y$ there is equal to the reduction of $-t'$.

\smallskip

\underline{Case 2: $\PP \nmid 2$.}

Here we cannot have~$\PP | t'$ or $\PP | (t'-2)$ because $\frac{t'-2}{t'}$ is in $ \OO_F^\times$.

If $\PP | (t'^2+2t'-4)$ and $\PP \nmid (t'^2-3t'+1)$, then $\PP \nmid 5$.
By Lemma 3.1 of \cite{Ba}, there exist $\mu,\nu$ in $\OO_F$ such that $\mu c_4 + \nu \Delta = 2^4 5$.
It follows that~$v_\PP(c_4)=0$, so that the above Weierstra\ss\
equation of $E'$ is minimal at~$ \PP $.
The point $P=(0,0)$ does not reduce to a singular point.

If $\PP | (t'^2-3t'+1)$ and $\PP \nmid (t'^2+2t'-4)$, then the
Weierstra\ss\ equation of~$E'$ is minimal at~$ \PP $ as above.
The reduction of $2P = (t'(t'-1)(t'-2), -t'(t'-1)^2(t'-2)^2)$ is regular.

Finally, if $\PP | (t'^2-3t'+1)$ and $\PP | (t'^2+2t'-4)$, then $\PP | 5$ and $\PP | (t'+1)$.
By the transformation in \cite[Ch III, Section 1]{Si}, $E'$ is isomorphic
to the curve defined by
\begin{equation*}
y^2 = x^3 - 27c_4x - 54c_6
\end{equation*}
with $2P$ corresponding to $(3(t^6 - 4t^5 + 20t^3 - 28t^2 + 8t + 4), -108t(t-1)^2(t-2)^2)$.
Then one shows as when~$\PP | 7$ in~(1) that~$2P$ does not reduce to a singular point.
\end{proof}

We now describe all cubic fields~$ F $ generated by an element~$ u $ such that both
$ u $ and $ 1-u $ are in~$ \OO_F^\times $
(i.e., $ u $ is an \emph{exceptional unit} of $\OO_F$). Note
that this is the exact condition on~$ u = t $ in part~(1) of Theorem~\ref{thm:integrality}
and on $u = \frac{1}{t}-1$ in part~(2).
We also describe all quartic fields~$ F $ generated by an element~$ u $ such
that both $ u $ and $ \frac{u-1}{u+1} $ are in~$ \OO_F^\times $,
which applies to part~(3) for~$u = 1-2t$.

\begin{table}[t]
\caption{\label{table:quarticfields} The coefficients $b, c, \v$ and the Galois groups $\Gal$ of the splitting fields of $X^4+aX^3+bX^2+cX+\v$ for general $a$.}
\centering
\begin{tabular}{|c|c|c|c||c|c|c|c|}
\hline
$b$ & $c$ & $\v$ & $\Gal$ & $b$ & $c$ & $\v$ & $\Gal$\\
\hline
$-2$ & $-a \pm 1$ & $1$ & $S_4$ & $0$ & $-a\pm 1$ & $-1$ & $S_4$  \\
$-2$ & $-a \pm 2$ & $1$ & $S_4$ & $0$ & $-a\pm 2$ & $-1$ & $S_4$ \\
$-2$ & $-a \pm 4$ & $1$ & $S_4$ & $0$ & $-a\pm 4$ & $-1$ & $D_4$  \\
$-2$ & $-a \pm 8, 2|a$ & $1$ & $S_4$  & $0$ & $-a\pm 8$, $2|a$ & $-1$ & $S_4$ \\
$-2$ & $-a \pm 16, 4|a$ & $1$ & $S_4$ & $0$ & $-a\pm 16$, $a\equiv 2\pmod 4$ & $-1$ & $S_4$ \\
$-2\pm 1$ & $-a$ & $1$ & $D_4$ & $\pm 1$ & $-a$ & $-1$ & $S_4$ \\
$-2\pm 2$ & $-a$ & $1$ & $D_4$ & $\pm 2$ & $-a$ & $-1$ & $S_4$ \\
$2$ & $-a$ & $1$ & $C_2 \times C_2$ & $\pm 4$ & $-a$ & $-1$ & $S_4$ \\
$-6$ & $-a$ & $1$ & $C_4$ & & & & \\
$-2\pm 8$ & $-a, 2|a$ & $1$ & $D_4$ & $\pm 8$ & $-a$, $2|a$ & $-1$ & $S_4$ \\
$-2\pm 16$ & $-a, 4|a$ & $1$ & $D_4$ & $\pm 16$ & $-a$, $a\equiv 2\pmod 4$ & $-1$ & $S_4$ \\
\hline
\end{tabular}
\end{table}

\begin{lemma}\label{lemma:unit}
(1)
For every integer~$ a $, and all $ \v , \v'$ in $ \{\pm1\} $, the polynomial $ f_a(X) = X^3+a X^2- (a+\v + \v' + 1) X + \v $ is irreducible in~$\Q[X]$.
A cubic field $F$ has an element~$u$ such that $F=\Q(u)$ and both $u$ and $1-u$ are in $\OO_F^\times$ precisely when~$u$ is a root of some
$ f_a(X) $. Moreover,~$ F/\Q $ is cyclic if and only if~$ \v = \v' = 1 $ or~$ |  2 a - \v + \v' + 3 | = 7 $.

(2)
For each of the 40 families indexed by~$ b $, $ c $ and $ \v $ as listed in Table~\ref{table:quarticfields},
with~$ a $ an arbitrary integer unless a condition on it is given in the column for~$ c $,
the polynomial~$ f_a(X) = X^4+aX^3+bX^2+cX+\v$ is irreducible in~$ \Q[X] $,
except for the 28 polynomials that only have
$ X^2 +1 $, $ X^2 \pm X - 1 $, $ X^2 \pm 2 X - 1 $, or $ X^2 \pm 4 X - 1 $
as factors.
For each family, with the possibility of finitely many exceptions,
the Galois group of the splitting field of the
polynomial is as listed in the table.

A quartic field~$ F $ has an element $ u $ such that $ F = \Q(u) $
and both $ u $ and $ \frac{u-1}{u+1} $ are in~$ \OO_F^\times $
precisely when $ u $ is a root of one of the irreducible~$ f_a(X) $.

\smallskip

Moreover, for $ F= \Q(u) $ with $ u $ a root of an $ f_a(X) $
in any of the families in~(1) and~(2),
the number of roots of unity in $ F $ is 2, 4, 6, 8 or~12,
and~$ F $ is totally real for $ \abs{a} \gg 0 $.
\end{lemma}

\begin{proof}
(1)
If both~$ u $ and $ 1-u $, or, equivalently, $ u - 1 $, are in $ \OO_F^\times $ with $ F $
a cubic field then $ u $ is a root of a polynomial~$f(X) = X^3+aX^2+bX+\v$ for integers~$ a $
and $ b $ with $\v =\pm 1$, and $ u - 1 $ of~$ g(X) = f(X+1) $.
Because $ u \ne \pm 1 $, $ F = \Q(u) $, so $ f(X) $ is irreducible in~$ \Q[X] $, hence the
same holds for $ g(X) $. With $ u - 1 $ in~$ \OO_F^\times $,
it follows that $ a + b + \v = f(1) = g(0) = - \v' $ for some $ \v' = \pm 1 $,
so that $ f(X) $ is of the required shape.
Conversely, from reducing its coefficients modulo~2 one sees any such $ f(X) $ is irreducible,
and the conditions on its coefficients imply that, for a root $ u $ of $ f(X) $,
we have~$ F = \Q(u) $ of degree 3 over~$ \Q $ with~$ u $ and $ u-1 $ in~$ \OO_F^\times $.

The discriminant of $ f(X) $ is~$ M^2 + 8 ((\v + 1) (\v' + 1) - 4) $
with
\begin{equation*}
 M = a^2+(-\v+\v'+3) a + (\v + 2) (2 \v' + 1 )
\,.
\end{equation*}
So this discriminant equals~$ M^2 $ if $ \v = \v' = 1 $. In the other cases
it equals~$ M^2 - 32 $, which is not a square whenever~$ M > 16 $,
and for~$ M \le 16 $ it can be checked directly.

If~$ F = \Q(u) $ with~$ u $ a root of any of the $ f_a(X) $,
then $ F $ is cubic, hence has only two roots of unity. From
the discriminant it follows that $ F $ is totally real for~$ \abs{a} \gg 0 $.
For later use we observe that, in fact,
the holomorphic implicit function theorem implies that for $|a| \gg 0$,
the roots of the polynomial $ a^{-1} f_a(X) $ are given by
$ u_0 = \v a^{-1} - \v (\v' + 1) a^{-2} + \dots $,
$ u_1 = 1 + \v' a^{-1} +  (\v-2) \v' a^{-2} + \dots $, and hence also $ u_\infty = -a - u_0 - u_1 $,
giving the three embeddings of $ F $.
Note that $ \log\abs{u_0} \sim - \log\abs{a} $,
$ \log\abs{u_1-1} \sim - \log\abs{a} $
and $ \log\abs{u_\infty} \sim \log\abs{a} $.

(2)
If both~$ u $ and~$ \frac{u-1}{u+1} $ are units in $ \OO_F $
with $ F $ a quartic field, then~$ u $ is a root of some $ f(X) = X^4 + a X^3 + b X^2 + c X + \v $ for
integers~$ a, b $ and $ c $, with~$ \v = \pm1 $. Then
$ \frac{u-1}{u+1} $ is a root of~$ g(X) = (1-X)^4 f(\frac{1+X}{1-X}) $, which
has leading coefficient $ -a + b -c + \v + 1 $ and constant
term $ a + b + c + \v + 1 $. If $ f(X) $, hence $ g(X) $, is irreducible,
then $ \frac{u-1}{u+1} $ is a unit of $\OO_F$ if and
only if $ g(X) $ scales to a monic polynomial in $ \Z[X] $ with constant term~$ \pm 1 $.
Comparing leading coefficient and constant term, we see that
either~$ a + c = 0 $, in which case $ g(X) $ equals
\begin{equation*}
 ( b + \v + 1) X^4 + 4 (-a - \v + 1 ) X^3 + 2( - b + 3 \v + 3) X^2 + 4 ( a -  \v + 1) X + (b + \v + 1)
,
\end{equation*}
or~$ b + \v + 1 = 0 $, so that it equals
\begin{equation*}
 (-a - c) X^4 + 2 (- a + c - 2 \v + 2 ) X^3 + 8( \v + 1) X^2 + 2 (a - c - 2 \v + 2) X + (a + c)
.
\end{equation*}
Demanding that all coefficients are divisible by the leading
one now gives the listed~40 families of polynomials $ f(X) $ with the desired scaling
of~$ g(X) $.

If $ f(X) $ is reducible then it can satisfy the scaling condition
imposed of the associated $ g(X) $ in only~28 cases.
To see this, note that~$ f(X) $ cannot have a factor~$ X \pm 1 $ as it would
lead to~$ g(X) $ being of degree~3 or having constant
term~0. If $ f(X) = q_1(X) q_2(X) $ with the~$ q_i(X) $ in~$ \Z[X] $
monic, quadratic, and with constant term $ \pm1 $, then
each~$ (1-X)^2 q_i(\frac{1+X}{1-X}) $ must scale
to such a quadratic. It is easy to check this holds for precisely the 7 quadratic
factors as listed in the lemma.

We sketch how the Galois group for each of the 40 families was determined using~\cite[Corollary~4.3]{Co}.
The discriminant of $ f = f_a $ is in $ \Z[a] $ of degree 6
with leading coefficient 4. So if it is not a square in $ \Z[a] $
then its value is not a square for all $ |a| \gg 0 $. It is a
square in~$ \Z[a] $ only when $ b =2 $ and $ \v = 1 $.

The cubic resultant $ R_3(Y) $ of $ f_a $ in $ \Z[Y,a] $ defines an affine
part of an elliptic curve for 28 families (not including the
case $ b = 2 $ and $ \v = 1 $). So, in each of those families, it has integer
zeroes for only finitely many~$ a $ by Siegel's theorem,
and for all other~$ a $ the Galois group is $ S_4 $.
For the remaining 12 families it splits into three linear factors
in~$ \Z[Y,a] $ only for~$ b = 2 $ and $ \v = 1 $, i.e., the case of~$ C_2\times C_2 $.
For the other 11 cases it is the product of a linear factor and
of a quadratic~$ Y^2 + B Y - a^2 + C $ for integers~$ B $
and~$ C $, or~$ Y^2 + ( 2 \pm  a) Y + 8  $. In all cases
this quadratic remains irreducible for almost all $ a $ in $ \Z $ as
its discriminant is~$ 4 a^2 + B^2 - 4 C $ or $ a^2 \pm 4 a - 28 $.

The remaining checks in those 11 cases to determine if the Galois group
is $ C_4 $ or~$ D_4 $ come down to verifying if an element of~$ \Z[a] $
of degree~8 with leading coefficient~1 or~4 evaluates to a square for an integer $ a $.
This is again the case infinitely often only when the expression
is a square in $ \Z[a] $ already, which happens for both checks
involved only for~$ b = -6 $ and $ \v = 1 $.

Because $ a^{-1}f_a(X) $ is of the form $ a^{-1} X^4 + X^3 + a^{-1} b X^2 + (-1 + a^{-1} c' ) X + a^{-1} \v $,
for~$|a| \gg 0$ it has real roots
$ u_{-1} = -1 - \frac12 (b-c'+\v+1) a^{-1} + \dots $,
$ u_0 = \v a^{-1} + \v c' a^{-2} + \dots $
and
$ u_1 = 1 - \frac12 (b + c' + \v + 1) a^{-1} + \dots $,
and together with the root~$u_\infty = -a - u_{-1} - u_0 - u_1 $
this gives four real embeddings of $ F $.
For later use we note that the coefficient of $a^{-1} $ in~$ u_{-1} $, $ u_0 $ and~$ u_1 $ is non-zero
for each of the 40 families of Table~\ref{table:quarticfields}.

Using the criteria for the roots of a quartic equation~\cite{R22},
for each family we can find an explicit lower bound on $ \abs{a} $ for
which all roots of $ f_a(X) $ are real.
Direct calculation for the remaining irreducible $ f_a(X) $ then gives the following exceptional cases in which the number of roots of unity $ \omega $ in~$F$ is greater than 2:
\begin{alignat*}{3}
&X^4+aX^3-aX+1,   \quad    && a=\pm 2\ (\omega = 4),\ 0 \ (\omega = 8); \\
&X^4+aX^3+2X^2-aX+1, \quad && a=\pm 3 \ (\omega = 6),\ \pm 2 \ (\omega = 12),\ \pm 1 \ (\omega = 6); \\
&X^4+aX^3+6X^2-aX+1, \quad && a=\pm 4 \ (\omega = 4),\ 0 \ (\omega = 8); \\
&X^4+aX^3-X^2-aX+1, \quad  && a=\pm1 \ (\omega = 6),\ 0 \ (\omega = 12).
\end{alignat*}
The cases where the polynomial is reducible are easily done directly.
\end{proof}

\begin{remark}
It is interesting to note that the family $X^3 + aX^2 - (a+3)X + 1$ in (1) defines the simplest cubic fields as in \cite{Sh},
and the family $X^4 + aX^3 - 6X^2 - aX + 1 $ in~(2),
which is irreducible precisely when $ a\neq 0, \pm3$,
defines the simplest quartic fields as in \cite{Gr}. The class groups and unit groups of these fields are well studied~\cite{Sh, Gr,Wa,L}.
Mahler measures of polynomials over these fields will be considered in~\cite{BLW}.
\end{remark}

\begin{remark}\label{remark:infinitefields}
For a given integer $a$, Hoshi and Miyake \cite[Corollary 5.6]{HM09} proved that,
for given $ a $ in~$ \Z $, there are only finitely many $b$ in $\Z$ such that the splitting fields of $X^3 + aX^2 - (a+3)X + 1$ and of~$X^3 + bX^2 - (b+3)X + 1$ coincide.
So in this family we have infinitely many different fields. One can
prove similar results for the other families of fields in Lemma \ref{lemma:unit} using the method of loc.~cit.
\end{remark}

\begin{remark} \label{remark:isofields}
In Lemma~\ref{lemma:unit}(1), if $ u $ and $ 1-u $ are in~$ \OO_F^\times $, then
$ v $ and $ 1-v $ are in~$ \OO_F^\times $ for~$ v $ in the orbit of~$ u $ under
the dihedral subgroup of order~6 of~$ \mathrm{PGL}_2(\Z) $ generated
by~$ (\begin{smallmatrix} -1 & 1 \\ 0 & 1 \end{smallmatrix}) $ and $ ( \begin{smallmatrix} 0 & 1 \\ 1 & 0 \end{smallmatrix} ) $.
So $ v $ equals~$ u^{\pm1} $, $ 1 - u^{\pm1} $ or~$ (1-u^{\pm1})^{-1} $,
and~$ F = \Q(v) $. Therefore we get isomorphic~$ F $ from the polynomials
in several ways.

If~$ \v = \v' = 1 $, then~$ f_a(X) = X^3 + a X^2 - (a+3) X + 1 $
satisfies
\begin{equation*}
f_a(X) = - X^3 f_a(1-X^{-1}) = (X-1)^3 f_a(1/(1-X))
,
\end{equation*}
and the maps~$ u \mapsto 1-u^{-1} $
and $ u \mapsto (1-u)^{-1} $ give the non-trivial elements in the Galois group of $ F $.
We also have $ X^3 f_a(X) = f_{-a-3}(X) $, so all fields are obtained for~$ a \ge -1 $.
Denoting the remaining three families of polynomials by~$ f_a(X) $
(again), $ g_a(X) $ and~$ h_a(X) $, we have
\begin{alignat*}{1}
f_a(X) & = X^3 + a X^2 - (a+1) X + 1 \\
g_a(X) & = X^3 + a X^2 - (a+1) X - 1 = - f_{-a-3}(1-X) \\
h_a(X) & = X^3 + a X^2 - (a-1) X - 1 = X^3 f_{-a-2}(1-X^{-1})
\end{alignat*}
as well as the identity~$ f_a(X) =  X^3  f_{-a-1}(X^{-1}) $, so all non-cyclic fields are obtained from~$ f_a(X) $ for~$ a \ge 0 $.
But if we let $ t $ be obtained from~$ u $ by the formulae
just before Lemma~\ref{lemma:unit} and consider the resulting
curves in Theorem~\ref{thm:integrality} (as we shall do in Theorem~\ref{thm:independence}
below) then these identities do not always give corresponding
symmetries for the curves.

For $ N = 7 $ we have $ t = u $ and $ j(u)  = j (1-u^{-1}) = j ((1-u)^{-1}) $,
the other three possibilities for~$ v $ all giving the same different rational
function in~$ u $.
For~$ N = 8 $ we have~$ t = (u+1)^{-1} $, and viewing the $ j $-invariant
as a function of~$ u $, the six~$ v $ pair up to give three different
rational functions in~$ u $. One has~$ j(u) = j(u^{-1}) $,
and transforming $ u $ into $ u^{-1} $ corresponds to transforming~$ t $ into~$ 1-t $.
So these identities match the discussion for~$ N = 7 $ and~8 in Remark~\ref{remark:isoN} .
For $ N = 10 $ we have that if $ u $ and~$ \frac{u-1}{u+1} $
are both in~$ \OO_F^\times $, then the same applies to~$ v $
and $ \frac{v-1}{v+1} $ if we take~$ v $ equal to~$ \pm u^{\pm1} $ or~$ \pm \bigl( \frac{u-1}{u+1} \bigr)^{\pm1} $.
These~$ v $ form the orbit of~$ u $ under the action of the dihedral subgroup of order~8
of~$ \mathrm{PGL}_2(\Q) $ generated
by~$ (\begin{smallmatrix} 1 & -1 \\ 1 & 1 \end{smallmatrix}) $
and $ ( \begin{smallmatrix} -1 & 0 \\ 0 & 1 \end{smallmatrix} ) $.
With~$ t = \frac{1-u}2 $ we again view the $ j $-invariant as a function~$ j(u) $
of~$ u $. Letting the dihedral group act, we find that the $ j(v) $ pair up to
give four different rational functions of~$ u $.
The identities include~$ j(u) = j(-u^{-1}) $, and because mapping $ u $ to $ -u^{-1} $
corresponds to mapping~$ t $ to~$ \frac{t-1}{2t-1}$, this again
matches Remark~\ref{remark:isoN}.

We did not fully investigate the identifications among the 40 families of fields listed in Table~\ref{table:quarticfields}
(let alone the identifications among the resulting families of curves)
under this dihedral group. But every such family
is mapped to itself by at least one element of order 2, limiting
the number of identifications.
The family of a polynomial $ f(X) = X^4 + a X^3 + b X^2 + c X + \v $ is determined by $ a+c $,
$ b $ and $ \v $. If $ \v = 1 $ then $ X^4 f(X^{-1}) $ is in
the same family, and if $ \v = -1 $ then this holds for~$ f(-X) $
if $ a+c = 0 $, or $ - X^4 f(-X^{-1}) $ if $ b = 0 $.
\end{remark}

\begin{remark}
We use~$ N = 7 $, 8 and~10 because then we can easily give $ F $-rational
points on the modular curve~$ X_1(N) $, which has genus~0.
For smaller~$ N $ the construction of Section~\ref{blochsection}
gives at most two linearly independent elements.
For the remaining cases~$ N =9 $ or~12 where~$ X_1(N) $ has genus~0,
the situation is a little different from the cases we treat
in this paper.
\end{remark}

\subsection{The main result} \label{subsec main result}
We can now formulate our main result for our families on the elements~$ S_{P,s} $
of Section~\ref{blochsection}.
(See Remark~\ref{remark:ST} for a variation involving other elements.)
In particular, in the next theorem,
for~$ |a| \gg 0 $, we have constructed as many linearly independent
elements in $ \INTE $ for the resulting curve~$ E $ as predicted by Beilinson's conjecture.
Note that by what we saw in Section~\ref{blochsection}, and since $ K_2 $ of a number field is a torsion group, we only need to consider the~$ S_{P,s} $ with $1 \le 2s \le N-1$.
We recall from Lemma~\ref{lemma:unit} that the field $ F $ is non-Abelian for~$ |a| \gg 0 $ in 3
of the cubic families and 38 of the quartic families.

\begin{theorem}\label{thm:independence}
Consider the following fields $ F $ with element~$ t $, parametrised by an integer $ a $,
and integers~$ N $.
\begin{enumerate}
\item
For fixed $\v, \v'$ in $\{\pm 1\}$, let $ f_a(X) = X^3+a X^2- (a+\v + \v' + 1) X + \v $ be as in Lemma~\ref{lemma:unit}(1),
let $ u $ be a root of $ f_a(X) $, and let $ F = \Q(u) $.
Then put $ t = u $ and $ N  = 7 $, or $ t = 1/(u+1) $ and $ N = 8 $.

\item
Let $ f_a(X) $ be in one of the families in Lemma~\ref{lemma:unit}(2),
as listed in Table~\ref{table:quarticfields},
let $ u $ be a root of $ f_a(X) $, and let $ F = \Q(u) $.
Then put $t=\frac{1-u}{2}$ and $ N  = 10 $.
\end{enumerate}
For such~$ F $, $ t $ and $ N $, if~$\Delta_N(t)$ in~Table~\ref{tab: tnf}
is defined and non-zero, let~$E/F$ be the elliptic curve
defined by the corresponding Tate normal form~\eqref{eqn:TNF}
with~$ f = f_N(t) $ and~$ g = g_N(t) $.
Then, with~$ P = (0,0) $:
\begin{itemize}
\item
the elements $\gcd(2,N) \cdot S_{P,s}$ for $ s = 1,\dots,N - 1 $ are in~$ \INTE $;

\item
the~$S_{P,s}$ with $s=1,\dots,\lfloor \frac{N-1}{2}\rfloor$ are linearly independent for~$\abs{a}\gg 0$,
and in fact, their Beilinson regulator $R=R(a)$ satisfies
\begin{equation*}
\qquad\quad
 \lim_{\abs{a}\rightarrow \infty}\frac{R(a)}{\log^{\lfloor \frac{N-1}{2}\rfloor} \abs{a}} =
 C_N \cdot \left|\det\left(\frac{N^4}{3}B_3\left(\left\{\frac{st}{N}\right\}\right)_{1 \le s,t \le \lfloor \frac{N-1}{2}\rfloor}\right)\right|,
\end{equation*}
where $B_3(X) = X^3-\frac{3}{2}X^2+\frac{1}{2}X$ is the third Bernoulli polynomial,
$\{x\}$ denotes the fractional part of $x \in \R$, and
$C_7=1$, $ C_8 = C_{10} = 4 $.
\end{itemize}
\end{theorem}

\begin{proof}
The statement about integrality follows from Theorem~\ref{thm:integrality}
and the result on the number of roots of unity in Lemma~\ref{lemma:unit}.
The proof of the statement on the limit behaviour of the Beilinson regulator
is involved and will be given in Section~\ref{section:independence}.
It implies the claimed linear independence because the above
determinant is non-zero.
\end{proof}

\begin{remark}
(1)
In the theorem, the condition that $ \Delta_N(t) $ is defined
and non-zero fails for precisely the following~$ N $ and~$ u $
defined in~(1) and~(2) of the theorem.
For~$ N = 7 $, if $ u $ is a root of $ X^3 - 8 X^2 + 5 X + 1 $.
For~$ N = 10 $, if $ u $ is a root of~$ X^2 + 4 X - 1 $ or of~ $ X^2 - X - 1 $.

(2)
We allow reducible~$ f_a(X) $ in~(2) in the theorem because, according
to the Beilinson conjectures, on the elliptic curve over the resulting
quadratic field~$ F $ the elements~$S_{P,1} $ through~$ S_{P,4} $
should satisfy linear dependencies. Verifying this (which we
did not try) may provide
interesting evidence for the conjectures.
\end{remark}

\begin{remark}
Lecacheux has constructed units in number fields using modular
curves~\cite{Le1,Le2,Le3}. Consider the Abelian covering $X_1(N)\rightarrow X_0(N)$.
For several values of $N$ such that $X_0(N)$ has genus $0$, she constructed functions~$ t $ on~$X_1(N)$ such that $ t $
and $ t^{-1} $ are both integral over~$\Z[a]$,
where $a$ is a certain Hauptmodul for $\Gamma_0(N)$.
Specialising $a$ to suitable integers, we get number fields~$F$,
each with a distinguished unit of~$ \OO_F $.
It would be interesting to study the group $K_2^T$ of elliptic curves with
points of order~13, 16 or 25 over the fields constructed by Lecacheux
in those papers.
\end{remark}

\section{Linear independence of the elements} \label{section:independence}
In this section, we prove the linear independence of the elements
$ S_{P,s} $ as stated in Theorem~\ref{thm:independence}. The proof
is analytical, and requires formulae for the regulator pairing~\eqref{regulator-pairing}
in terms of the elliptic dilogarithm, which we recall in Section \ref{reg elliptic dilog}.
We note that historically, these formulae were important in the development
of the conjectures on special values of $L$-functions.
We also need a careful analysis of the action of complex conjugation
on the fibres of the universal elliptic curve $\mathcal{E}_1(N)$, which we do
in Section \ref{cc E1N}.
The proof of linear independence
is carried out in Section \ref{end proof}, by finding the
limit behaviour of the regulator pairing~\eqref{regulator-pairing}
associated to the elements $S_{P,s}$ when $ \abs{a} $ goes to infinity.

\subsection{The elliptic dilogarithm} \label{reg elliptic dilog}
Let $E$ be an elliptic curve defined over $\C$,
and choose an isomorphism $E(\C) \cong \C/(\Z + \tau\Z)$,
where~$\tau$ is in the upper half-plane.
The exponential map $u \in \C \mapsto \mathrm{e}(u) = \exp(2\pi iu)$
identifies $ \C/(\Z + \tau\Z) $ with $ \C^\times/q^{\Z} $, where $ q = e^{2 \pi i \tau} $.

The Bloch-Wigner dilogarithm $D\colon\mathbb{P}^1(\C) \to \R$
is the unique continuous function satisfying
$D(z) = \im(\sum_{n = 1}^\infty z^n/n^2) + \operatorname{Arg}(1-z) \log\abs{z}$
for $ 0 < |z| \leqslant 1 $, $z \neq 1$, and $D(1/z) = -D(z)$ for every $z$.
It satisfies $ D( \bar z ) = - D(z) $.
Bloch's elliptic dilogarithm $ D_q $ is then defined by averaging $D$ over~$q^\Z$
as
\begin{equation*}
D_q\colon\C^\times/q^\Z \to \R, \qquad z \mapsto \sum_{n\in\Z} D(zq^n).
\end{equation*}
Similarly, we define $J(z) = \log\abs{z}\log\abs{1-z}$
and $J_q\colon\C^\times/q^{\Z} \to \R$ by
\begin{equation*}
J_q(z) = \sum_{n=0}^{\infty}J(zq^n) - \sum_{n=1}^{\infty}J(z^{-1}q^n)+\frac{1}{3}\log^2\abs{q}B_3\left(\frac{\log\abs{z}}{\log\abs{q}}\right),
\end{equation*}
where $B_3(X)$ is the third Bernoulli polynomial.
Using that $ B_3(X+1) = - B_3(-X) = B_3(X) + 3 X^2 $, one easily verifies
that $ J_q(q z) = J_q(z) $ and $ J_q(z^{-1}) = - J_q(z) $.
Then we define the function $R_q\colon\C^\times / q^\Z \to \C$ as
$ R_q = D_q - i J_q $.
From the above we have that $ R_q(z^{-1}) = - R_q(z) $ and $ \overline{R_q(z)} = - R_{\bar q}(\bar{z}) $.

It will be convenient to view the function $R_q$ as defined on the
additive group $\C / (\Z + \t \Z)$, by means of the isomorphism
$\mathrm{e}\colon\C / (\Z + \t \Z) \to \C^\times / q^\Z$ introduced at
the beginning of the section. We denote by $R_\t = D_\tau - i J_\tau\colon\C / (\Z + \t \Z) \to \C$
the resulting function. It is then clear that $ R_\t(-u) = - R_\t(u) $,
$ \overline{R_\tau(u)} = R_{-\overline{\tau}}(\overline{u}) $,
and $ R_{\t+n}(u) = R_\t(u) $ for any integer~$ n $.

The regulator pairing~\eqref{regulator-pairing}
can be computed using~$R_\tau$~\cite[Corollary~3.2]{LR},
as we shall now explain.
(In fact, this computation by Bloch in~\cite{Bl00} preceded Beilinson's
more general definition of the regulator map.)
Let $f$ and $g$ be non-zero rational functions on~$E$ with divisors
\begin{equation*}
(f) = \sum m_i(a_i), \quad (g) = \sum n_j(b_j).
\end{equation*}
We define the diamond operator by
\begin{equation*}
(f) \diamond (g) = \sum m_i n_j (a_i-b_j) \in \Z[E(\C)]^-,
\end{equation*}
where $\Z[E(\C)]^-$ denotes the quotient of $\Z[E(\C)]$ by the subgroup generated by $[P]+[-P]$ for all~$P$ in~$E(\C)$. The diamond operator induces a group homomorphism
\begin{equation*}
\C(E)^\times \otimes_{\Z} \C(E)^\times \rightarrow \Z[E(\C)]^-.
\end{equation*}
We extend $R_\t$ to a function on $\Z[E(\C)]$ by setting
$R_\t(\sum m_i (a_i)) = \sum m_i R_\t(a_i)$.
Since $ R_\t $ as an odd function on $ E(\C) $, it induces a map
$ R_\t\colon\Z[E(\C)]^- \rightarrow \C $.

Now let $\gamma_0$
be the path from $0$ to $1$ in $E(\C) \cong \C/(\Z+\tau\Z)$.
Then for any $\gamma \in H_1(E(\C), \Z)$ and any non-zero
holomorphic 1-form $ \omega $ on $E$,
with associated periods~$\Omega_{\g_0} = \int_{\gamma_0} \omega$
and~$\Omega_\g = \int_{\g}\omega$, we have, if ~$\a = \sum_j \{f_j, g_j \} $,
\begin{equation}\label{eqn:reg}
\langle \g, \alpha \rangle = - \frac 1 {2\pi}\im\Bigl(\frac{\Omega_\g}{y_\tau \Omega_{\g_0}} \sum_j R_\tau((f_j) \diamond (g_j))\Bigr) \qquad (y_\tau = \im(\tau)).
\end{equation}
For fixed $ \g $, the quotient $ \Omega_\g / \Omega_{\g_0} $
is determined by the parametrisation. But it determines how $ D_\tau $ and $  - i J_\tau $
(which behave differently) contribute to~\eqref{eqn:reg}.

We shall use~\eqref{eqn:reg} for $ \a = S_{P,s} $ and $\gamma = \gamma_0$.
From the description of $S_{P,s}$ we find the element
\begin{equation} \label{Spdiamond}
\begin{aligned}
( N (s P) - N (O)) \diamond \Bigl( \sum_{t=0}^{N-1} N(tP) - N^2 (O) \Bigr)
& =
 ( N (s P) - N (O)) \diamond ( - N^2 (O) )
\\
& =
 - N^3 (s P) + N^3 (O)
\end{aligned}
\end{equation}
in $ \Z[E(\C)]^- $, so that
\begin{equation}\label{eqn:regSP}
\langle \gamma_0,  S_{P,s} \rangle
= \frac{N^3}{2\pi y_\tau} \im(R_\tau(sP))
= - \frac{N^3}{2\pi y_\tau} J_\tau(sP)
.
\end{equation}

Finally, we recall the Fourier expansion of $R_\tau(u)$ with respect to $\tau$ from~\cite[Section~3]{DI}
as we shall use it in the proof of Theorem~\ref{thm:independence}.
Assume that $u=a+b\tau$ with $0 \leqslant a, b < 1$,
and write $\tau = x+iy$ with $x \in \R$ and $y>0$. Then
\begin{align}
\begin{split}\label{eqn:Fourier}
D_\tau(u) &= -\frac{i}{2}\sum_{\substack{m-b, n \in \Z \\ n \neq 0}} \mathrm{e}(n a) \mathrm{e}(m n x) e^{-2 \pi|m n| y}\left(2 \pi y \frac{n|m|}{|n|^2}+\frac{n}{|n|^3}\right),\\
J_\tau(u) &= \frac{4 \pi^2 y^2}{3} B_3(b) - \pi y \sum_{\substack{m, n+b \in \Z \\ m, n \neq 0}} \mathrm{e}(-m a) \mathrm{e}(m n x) \frac{n}{|m|} e^{-2 \pi|m n| y}.
\end{split}
\end{align}

\subsection{Complex conjugation on the universal elliptic curve} \label{cc E1N}
Recall from Section \ref{subsection:parametrization} that the
isomorphism $\nu' \colon \Gamma_1(N) \backslash \H \to Y_1(N)(\C)$ is
compatible with complex conjugation, where the complex conjugation
on $\Gamma_1(N) \backslash \H$ is induced by $c(\tau) = -\bar{\tau}$.
Moreover, in the notations of Section \ref{subsection:parametrization},
the complex conjugation on $\mathcal{E}_1(N)(\C)$ sends $[\tau, z]$
to $[-\bar\tau, \bar z]$.

We now describe the complex conjugation on the fibres of $\mathcal{E}_1(N)(\C)$.

\begin{lemma} \label{lem:cc}
Let $\tau$ in $\H$ be such that $\nu'(\tau) \in Y_1(N)(\R)$.
Write $c(\tau) = A \tau$ with $A =
(\begin{smallmatrix} s & t \\ u & v \end{smallmatrix}) \in \Gamma_1(N)$.
Identifying the fibre of $\mathcal{E}_1(N)(\C)$ over $\nu'(\tau)$ with
$\C/(\Z+N\tau\Z)$ as in Section \ref{subsection:parametrization},
the complex conjugation on this fibre
is given by $z \mapsto -(u\tau+v) \bar z$.
\end{lemma}

\begin{proof}
Note that $c(W_N(\tau)) = W_N(A \tau) = A' W_N(\tau)$
with $A' = (\begin{smallmatrix} v & -u/N \\ -Nt & s \end{smallmatrix})$
in $\Gamma_1(N)$. Let $z \in \C/(\Z+N\tau\Z)$, corresponding to
the point $[W_N(\tau), \frac{z}{N\tau}]$ in $\mathcal{E}_1(N)(\C)$. Its complex conjugate is
\begin{equation*}
\Bigl[A' W_N(\tau), \frac{\bar{z}}{N\bar\tau}\Bigr] = \Bigl[W_N(\tau), \Bigl(\frac{t}{\tau}+s\Bigr) \frac{\bar{z}}{N\bar\tau}\Bigr],
\end{equation*}
corresponding to the point $(s\tau+t) \frac{\bar{z}}{\bar\tau}$ in $\C/(\Z+N\tau\Z)$.
We conclude by noting that $c(\tau) = (s\tau+t)/(u\tau+v)$.
\end{proof}

We shall need the following lemma about the real points of $X_1(N)$.
For a complete description of $X_1(N)(\R)$, see \cite{Sn}.
For any $\alpha, \beta$ in $\Q \cup \{i\infty\}$,
$\alpha \neq \beta$, we denote by $\{\alpha, \beta\}$ the
open geodesic from $\alpha$ to $\beta$ in $\H$.

\begin{lemma} \label{lem X1NR}
Let $N \geqslant 4$, and $\alpha \neq \beta$ in $\Q \cup \{i\infty\}$
such that $\{\alpha, \beta\}$ is contained in
the real locus of $\Gamma_1(N) \backslash \H$.
Choose a matrix $B \in \mathrm{SL}_2(\Z)$ such that
$N\alpha = B \cdot \infty$, and write $N\tau = B \cdot \mu$ with $\mu \in \H$.
Then for $\tau \in \{\alpha, \beta\}$, we have  $\re(\mu) \in \frac12 \Z$.
\end{lemma}

\begin{proof}
By \cite[Theorem 3.1.1]{Sn}, there exists $A = (\begin{smallmatrix}
s & t \\ u & v \end{smallmatrix})$ in $\Gamma_1(N)$ such that
$c(\tau) = A \tau$ for every $\tau \in \{\alpha, \beta\}$.
Then $c(N\tau) = A' (N\tau)$ with $A' = (\begin{smallmatrix}
s & Nt \\ u/N & v \end{smallmatrix}) \in \mathrm{SL}_2(\Z)$.
Let $C = (\begin{smallmatrix} 1 & 0 \\ 0 & -1 \end{smallmatrix})$.
For $\tau \in \{\alpha, \beta\}$, we have
\begin{equation*}
c(\mu) = c(B^{-1} \cdot N\tau) = CB^{-1}C \cdot c(N\tau) = B_0 \mu
\end{equation*}
with $ B_0 = C B^{-1} C A' B$ in~$ \mathrm{SL}_2(\Z) $.
Note that $\mu$ approaches $i \infty$
when $\tau$ approaches $\alpha$. Taking the limit when $\mu \to i\infty$,
we get $B_0 \cdot i \infty = i \infty$, so that $B_0$ must be of the form
$\pm (\begin{smallmatrix} 1 & k \\ 0 & 1 \end{smallmatrix})$ with $k \in \Z$.
This shows that $c(\mu) = \mu+k$, hence the result.
\end{proof}

Under the assumptions of Lemma \ref{lem X1NR}, for
$\tau \in \{\alpha, \beta\}$, there is an isomorphism
$\C/(\Z+N\tau\Z) \cong \C/(\Z+\mu\Z)$. We now make explicit
the action of complex conjugation on $\C/(\Z+\mu\Z)$
under this isomorphism.

\begin{lemma} \label{lem:cc2}
Let $N \geqslant 4$, and $\alpha, \beta$ in $\Q \cup \{i\infty\}$ with
$\beta \neq \alpha, i\infty$, such that the geodesic $\{\alpha, \beta\}$
is contained in the real locus of $\Gamma_1(N) \backslash \H$.
On this geodesic, write $c(\tau) = A\tau$ with $A = (\begin{smallmatrix}
s & t \\ u & v \end{smallmatrix}) \in \Gamma_1(N)$.
Under the isomorphism $\C/(\Z+N\tau\Z) \cong \C/(\Z+\mu\Z)$, the complex
conjugation on $\C/(\Z+\mu\Z)$ is given by $z \mapsto \varepsilon \bar{z}$,
where
\begin{equation*}
\varepsilon = \begin{cases} -1 & \textrm{if } \alpha = i\infty, \\
-1 & \textrm{if } \alpha \neq i\infty  \textrm{ and } (\b-\a)u >0, \\
\phantom{-}1 & \textrm{otherwise.}
\end{cases}
\end{equation*}
\end{lemma}

\begin{proof}
Keeping the notations of Lemma \ref{lem X1NR}, write $N\tau = B \cdot \mu$
with $B = (\begin{smallmatrix} a_0 & b_0 \\ c_0 & d_0 \end{smallmatrix})$.
We fix the isomorphism
\begin{equation} \label{eq:cc2}
\begin{split}
\phi_B\colon\C/(\Z+N\tau\Z) & \to \C/(\Z+\mu\Z) \\
z & \mapsto (c_0 \mu + d_0) z.
\end{split}
\end{equation}
Note that the complex conjugation on $\C/(\Z+\mu\Z)$ is an antiholomorphic
isomorphism fixing $0$, hence must be of the form $z \mapsto \pm \bar z$
for $\im(\mu) \gg 0$. By continuity, this is true for any $\mu$, and the sign
is independent of $\tau \in \{\alpha, \beta\}$. Using Lemma \ref{lem:cc} and
\eqref{eq:cc2}, we compute that the complex conjugation on $\C/(\Z+\mu\Z)$
is given by
\begin{equation} \label{eq:cc3}
z \mapsto - \Bigl(\frac{u}{N} \frac{a_0 \mu + b_0}{c_0 \mu + d_0} + v\Bigr) \frac{c_0 \mu + d_0}{c_0 \bar\mu + d_0} \bar z.
\end{equation}
To determine the sign, it suffices to let $\mu \to i\infty$. Noting that
$\mu/\bar\mu \to -1$, we obtain in the case $c_0 \neq 0$
(equivalently, $\alpha \neq i\infty$):
\begin{equation*}
\varepsilon = \frac{u a_0}{N c_0} + v = \alpha u + v.
\end{equation*}
Furthermore, the equation $c(\tau) = A\tau$ can be solved explicitly:
for $s \neq v$ there is no solution, while for $s=v$, this locus is a
semi-circle passing through the points $(-s \pm 1)/u$. It follows that
$\alpha u + v = (-s \pm 1 + v) = \pm 1$, and the sign can be
determined according to the relative positions of $\alpha, \beta$
and the sign of $u$. Finally, in the case $\alpha = i \infty$, we have
$B = I_2$, and the equation $c(\tau)=A\tau$ implies $u=0$ and $v=1$,
so that \eqref{eq:cc3} simplifies to $z \mapsto -\bar z$.
\end{proof}

\subsection{End of the proof of Theorem~\ref{thm:independence}} \label{end proof}
In all cases, we shall determine the limit behaviour of
the Beilinson regulator in~\eqref{eqn:L-reg}
for the elements $S_{P,s}$ as~$ \abs{a} \to \infty $.
This implies the desired linear independence.

Let $E_t$ denote the elliptic curve associated to the parameter
$a \in \Z$, defined over the corresponding number field $F = \Q(t)$.
We make explicit the Beilinson regulator~$ R = R(a) $ in~\eqref{eqn:L-reg}
for the elements~$ \gcd(2,N) \cdot S_{P,j} $ with $ j = 1, \dots, \lfloor \frac{N-1}{2}\rfloor $.
As stated in Lemma~\ref{lemma:unit},
the field~$ F $ is totally real for~$ \abs{a} \gg 0 $, and its
embeddings correspond to the (real) roots $u$ of the polynomial~$ f_a(X) $
defining~$ F $. The asymptotics of these roots as $ \abs{a} \to \infty $
are also given in the proof of that lemma.
We shall indicate which embedding is considered by adding the limit of the root~$ u $
as subscript.

The required homology group $H_t = H_1(E_t(\C), \Z)^-$
in Section~\ref{conjecture} is then the direct sum
of the groups $H_i = H_1(E_i(\C), \Z)^- $, with~$ i $ running
over the limits  of~$ u $.
We find an explicit generator~$ \g_i $ of $H_i$ by combining the tools from Section~\ref{cc E1N} with an explicit description of $X_1(N)(\R)$
due to Snowden \cite{Sn}.
The Beilinson regulator $R(a)$ equals~$ \gcd(2,N) ^{[F:\Q]} \cdot |\! \det ( \langle \g_i, S_{P,j}  \rangle_X )_{i,j} |$,
where $ i $ runs through the values of the limits of~$ u $, and~$ j = 1, \dots, \lfloor \frac{N-1}{2}\rfloor $,
and all entries in the determinant depend on~$ a $.
As explained in Section~\ref{conjecture},
these entries are computed
by pulling back~$S_{P,s}$ to $ E_i $,
and pairing the result with~$ \g_i $ under the regulator pairing~\eqref{regulator-pairing} on $E_i$.
We emphasise that we compute the regulator with respect
to a $\Z$-basis of $H_t$, not just a $\Q$-basis of $H_t \otimes \Q$.
In all cases, the condition $|a| \to \infty$ corresponds to the modular
parameter $\tau$ in $\H$ approaching suitable cusps along the real locus of
$\Gamma_1(N) \backslash \H$. Using Sections~\ref{reg elliptic dilog}
and~\ref{cc E1N}, we determine the corresponding behaviour of the
function $R_\tau$, to get the desired limit behaviour of the regulator pairing.

(1)
We let $ u_0 $, $ u_1 $ and $ u_\infty $ for $ |a| \gg 0 $ be as in the proof of Lemma~\ref{lemma:unit}(1).

We begin with the case $ N = 7 $. Let $1 \leqslant s \leqslant 3$.
We denote by $ t_0 = u_0 $, $ t_1 = u_1 $ and $ t_\infty = u_\infty $
the corresponding embedded values of $t$. When $|a| \to \infty$,
we have $t_0 \to 0$, $t_1 \to 1$ and $t_\infty \to \infty$.
By Table \ref{table cusps N=7}, this corresponds to $\tau$
approaching the cusps $\frac27$, $i\infty$, $\frac37$ respectively,
along the real locus of $\Gamma_1(7) \backslash \H$.

In what follows, we use the notations of
Lemmas \ref{lem X1NR} and~\ref{lem:cc2}.

\mybullet $\tau \to i\infty$.
We have $B=I_2$, $7\tau = \mu$ and by Lemma \ref{lem:cc2},
the complex conjugation on $E_{t_1} \cong \C/(\Z+7\t\Z)
= \C/(\Z+\mu\Z)$ is given by
$z \mapsto -\bar z$. Therefore a generator $\gamma_1$ of
$H_1(E_{t_1}, \Z)^-$ corresponds to the loop from $0$ to $1$.
The regulator pairing of~$S_{P,s}$ and $\gamma_1$
is then computed using \eqref{eqn:regSP} with the parametrisation
$E_{t_1} \cong \C/(\Z+\mu\Z)$, under which the point $P$ corresponds
to $\mu/7$. We get
\begin{equation*}
\langle \gamma_1, S_{P,s} \rangle = - \frac{7^3}{2\pi y_\mu} J_\mu\Bigl(\frac{s \mu}{7}\Bigr).
\end{equation*}
The Fourier expansion~\eqref{eqn:Fourier} of $J_{\tau}$ provides the asymptotics
\begin{equation*}
J_\mu\bigl(\frac{s \mu}{7}\bigr) \sim_{|a| \to \infty} \frac{4\pi^2}{3} y_\mu^2 B_3\Bigl(\Bigl\{\frac{s}{7}\Bigr\}\Bigr).
\end{equation*}
Therefore
\begin{equation*}
\langle \gamma_1, S_{P,s} \rangle \sim_{|a| \to \infty} -2\pi \frac{7^3}{3} y_\mu B_3\Bigl(\Bigl\{\frac{s}{7}\Bigr\}\Bigr).
\end{equation*}

It remains to find the asymptotics of $y_\mu$ in terms of $|a|$.
By \eqref{t at oo}, we have the $q$-expansion $t \circ \nu' = 1 - q + O(q^2)$,
where $q = e^{2\pi i \tau}$, $\tau \to i\infty$. Recall also from
the proof of Lemma \ref{lemma:unit}(1) that $t_1 = 1 + \v' a^{-1} +  O(a^{-2})$.
Comparing these, we have $q \sim -\v' a^{-1}$, which implies
$y_\tau \sim \frac{\log|a|}{2\pi}$ and thus $y_\mu \sim \frac{7\log|a|}{2\pi}$. Finally, we have
\begin{equation} \label{eq regint 7-1}
\langle \gamma_1, S_{P,s} \rangle \sim_{|a| \to \infty} -\frac{7^4}{3} B_3\Bigl(\Bigl\{\frac{s}{7}\Bigr\}\Bigr) \log |a|.
\end{equation}

\mybullet $\tau \to \frac37$.
We write $7\tau = B\mu$ with $B = (\begin{smallmatrix} 3 & -1 \\ 1 & 0 \end{smallmatrix})$
and $\mu \to i\infty$. The isomorphism $\phi_B\colon\C/(\Z+7\tau\Z) \to \C/(\Z+\mu\Z)$
from \eqref{eq:cc2} sends $P=\tau$ to $\mu\tau = (3\mu-1)/7$.
We determine the complex conjugation on $\C/(\Z+\mu\Z)$ using
Lemma \ref{lem:cc2}. In the notations of this lemma, we have
$A = (\begin{smallmatrix} -13 & 6 \\ 28 & -13 \end{smallmatrix})$
on the geodesic $\{\frac37, \frac12\}$, and $A = (\begin{smallmatrix}
8 & -3 \\ -21 & 8 \end{smallmatrix})$ on $\{\frac37, \frac13\}$.
By Lemma \ref{lem:cc2}, on each side of the cusp $\frac37$, the complex
conjugation on $\C/(\Z+\mu\Z)$ is given by $z \mapsto -\bar z$, so that,
again, a generator $\gamma_\infty$ of $H_1(E_{t_\infty}, \Z)^-$
corresponds to the loop from $0$ to $1$.
As above, we get
\begin{equation*}
\langle \gamma_\infty, S_{P,s} \rangle = - \frac{7^3}{2\pi y_\mu} J_\mu\Bigl(\frac{s (3\mu-1)}{7}\Bigr).
\end{equation*}
which leads to
\begin{equation*}
\langle \gamma_\infty, S_{P,s} \rangle \sim_{|a| \to \infty} -2\pi \frac{7^3}{3} y_\mu B_3\Bigl(\Bigl\{\frac{3s}{7}\Bigr\}\Bigr).
\end{equation*}

We now estimate $y_\mu$ in terms of $|a|$. We use the modularity
property of the Weierstra\ss{} $\wp$-function for $\mathrm{SL}_2(\Z)$, which follows from the definition of $\wp$:
\begin{equation*}
\wp_{B_0\mu}\Bigl(\frac{z}{c_0 \mu + d_0}\Bigr) = (c_0 \mu + d_0)^2 \wp_\mu(z) \qquad \bigl(B_0 = (\begin{smallmatrix} a_0 & b_0 \\ c_0 & d_0 \end{smallmatrix}) \in \mathrm{SL}_2(\Z)\bigr).
\end{equation*}
Using this with $B_0 = B$ and $z = k\tau$, $1 \leqslant k \leqslant 6$, we obtain
\begin{equation*}
\wp_{7\tau}(k\tau) = \wp_{B\mu}\Bigl(\frac{k}{7} \cdot \frac{3\mu-1}{\mu}\Bigr) = \mu^2 \wp_\mu\Bigl(\frac{3k\mu-k}{7} \Bigr).
\end{equation*}
The Fourier expansion of this expression with respect to
$q^{1/7} = e^{2\pi i \mu/7}$ can be computed with \cite[Theorem 6.2(a)]{Sil94}.
Substituting into \eqref{eq t 7}, we get
\begin{equation*}
t \circ \nu'(\tau) = \frac{\wp_{7\tau}(2\tau)-\wp_{7\tau}(\tau)}{\wp_{7\tau}(3\tau)-\wp_{7\tau}(\tau)} = \zeta_7^5 q^{-1/7} + O_{\mu \to i\infty}(1)
\end{equation*}
with $\zeta_N := e^{2\pi i/N}$.
Combining this with the asymptotic~$t_\infty = -a + O(1)$ in the proof of Lemma \ref{lemma:unit}(1), we get
$y_\mu \sim \frac{7 \log |a|}{2\pi}$. Finally, we have
\begin{equation} \label{eq regint 7-2}
\langle \gamma_\infty, S_{P,s} \rangle \sim_{|a| \to \infty} -\frac{7^4}{3} B_3\Bigl(\Bigl\{\frac{3s}{7}\Bigr\}\Bigr) \log |a|.
\end{equation}

\mybullet $\tau \to \frac27$.
This case is similar to the case $\tau \to \frac37$, using the matrices
$B = (\begin{smallmatrix} 2 & -1 \\ 1 & 0 \end{smallmatrix})$,
$A = (\begin{smallmatrix} -13 & 4 \\ 42 & -13 \end{smallmatrix})$
on the geodesic $\{\frac27, \frac13\}$, and $A = (\begin{smallmatrix}
1 & 0 \\ -7 & 1 \end{smallmatrix})$ on $\{\frac27, 0\}$. This results in
\begin{equation} \label{eq regint 7-3}
\langle \gamma_0, S_{P,s} \rangle \sim_{|a| \to \infty} -\frac{7^4}{3} B_3\Bigl(\Bigl\{\frac{2s}{7}\Bigr\}\Bigr) \log |a|,
\end{equation}
where $\gamma_0$ denotes a generator of $H_1(E_{t_0}, \Z)^-$.

Combining \eqref{eq regint 7-1}, \eqref{eq regint 7-2} and \eqref{eq regint 7-3},
we obtain the desired limit behaviour
for the regulator of $S_{P,1}$, $S_{P,2}$ and $S_{P,3}$.
The linear independence of these elements for $|a| \gg 0$
boils down to the invertibility of the matrix
$\bigl(B_3(\{\frac{st}{7}\})\bigr)_{1 \le s,t \le 3}$,
which can be checked directly.

We now tackle the case $ N = 8 $, using the same method.
Let $1 \leqslant s \leqslant 3$.
This time $t = \frac{1}{u+1}$, and we denote by $t_0$, $t_1$ and
$t_\infty$ the values associated to $ u_0 $, $ u_1 $ and
$ u_\infty $. When $|a| \to \infty$, we have $t_0 \to 1$,
$t_1 \to \frac12$ and $t_\infty \to 0$. By Table \ref{table cusps N=8},
this corresponds to $\tau$ approaching the cusps
$i\infty$, $\frac14$, $\frac38$ respectively, along the real locus
of $\Gamma_1(8) \backslash \H$.

\mybullet $\tau \to i\infty$. This case is similar to the case $N=7$
and $\tau \to i\infty$. We get
\begin{equation} \label{eq regint 8-1}
\langle \gamma_0, S_{P,s} \rangle \sim_{|a| \to \infty} -\frac{8^4}{3} B_3\Bigl(\Bigl\{\frac{s}{8}\Bigr\}\Bigr) \log |a|,
\end{equation}
where $\gamma_0$ denotes a generator of $H_1(E_{t_0}, \Z)^-$.

\mybullet $\tau \to \frac38$. This case is similar to the case $N=7$
and $\tau \to \frac37$, using the matrices
$B = (\begin{smallmatrix} 3 & -1 \\ 1 & 0 \end{smallmatrix})$,
$A = (\begin{smallmatrix} -7 & 3 \\ 16 & -7 \end{smallmatrix})$
on the geodesic $\{\frac38, \frac12\}$, and $A = (\begin{smallmatrix}
17 & -6 \\ -48 & 17 \end{smallmatrix})$ on $\{\frac38, \frac13\}$.
By Lemma \ref{lem:cc2}, the complex conjugation on $\C/(\Z+\mu\Z)$
is given by $z \mapsto -\bar z$ in both cases, so a generator $\gamma_\infty$
of $H_1(E_{t_\infty}, \Z)^-$ is the loop from $0$ to $1$.
This leads to
\begin{equation} \label{eq regint 8-2}
\langle \gamma_\infty, S_{P,s} \rangle \sim_{|a| \to \infty} -\frac{8^4}{3} B_3\Bigl(\Bigl\{\frac{3s}{8}\Bigr\}\Bigr) \log |a|.
\end{equation}

\mybullet $\tau \to \frac14$.
This case is similar to the previous one, using the matrices
$B = (\begin{smallmatrix} 2 & -1 \\ 1 & 0 \end{smallmatrix})$,
$A = (\begin{smallmatrix} -7 & 2 \\ 24 & -7 \end{smallmatrix})$
on the geodesic $\{\frac14, \frac13\}$, and $A = (\begin{smallmatrix}
1 & 0 \\ -8 & 1 \end{smallmatrix})$ on $\{\frac14, 0\}$. This yields
\begin{equation} \label{eq regint 8-3}
\langle \gamma_1, S_{P,s} \rangle \sim_{|a| \to \infty} -\frac{8^4}{6} B_3\Bigl(\Bigl\{\frac{s}{4}\Bigr\}\Bigr) \log |a|,
\end{equation}
where $\gamma_1$ denotes a generator of $H_1(E_{t_1}, \Z)^-$.
Note that in this case $y_\mu \sim \frac{4\log|a|}{2\pi}$, while
$y_\mu \sim \frac{8\log|a|}{2\pi}$ in the other two cases.
This is due to the fact that the width of the cusp $\frac14$ is $2$,
while the other two cusps have width $1$.

From \eqref{eq regint 8-1}, \eqref{eq regint 8-2} and \eqref{eq regint 8-3},
we get the limit formula for the regulator of $2S_{P,1}$, $2S_{P,2}$ and $2S_{P,3}$.
The linear independence of these elements for $ \abs{a} \gg 0 $
follows from the fact that the determinant of the matrix
$\bigl(B_3(\{\frac{st}{8}\})\bigr)_{1 \le s,t \le 3}$ is non-zero.

(2)
We move on to the quartic case, corresponding to $N=10$.
Let $1 \leqslant s \leqslant 4$.
Now $t = \frac{1-u}2$, and we denote by $t_{-1}$, $t_0$, $t_1$ and
$t_\infty$ the values associated to $u_{-1}$, $ u_0 $, $ u_1 $ and
$ u_\infty $ as in the proof of Lemma~\ref{lemma:unit}(2) for $ |a| \gg 0 $.
When $|a| \to \infty$, we have $t_{-1} \to 1$,
$t_0 \to \frac12$, $t_1 \to 0$ and $t_\infty \to \infty$. By Table
\ref{table cusps N=10}, this corresponds to $\tau$ approaching the
cusps $i\infty$, $\frac15$, $\frac{3}{10}$ and $\frac25$ respectively,
along the real locus of $\Gamma_1(10) \backslash \H$.

\mybullet $\tau \to i\infty$.
As in the cases $N=7$ and $N=8$, we find
\begin{equation} \label{eq regint 10-1}
\langle \gamma_{-1}, S_{P,s} \rangle \sim_{|a| \to \infty} -2\pi \frac{10^3}{3} y_\mu B_3\Bigl(\Bigl\{\frac{s}{10}\Bigr\}\Bigr),
\end{equation}
where $\gamma_{-1}$ denotes a generator of $H_1(E_{t_{-1}}, \Z)^-$.

\mybullet $\tau \to \frac25$.
We have $B = (\begin{smallmatrix} 4 & -1 \\ 1 & 0 \end{smallmatrix})$
and $\phi_B\colon\C/(\Z+10\tau\Z) \to \C/(\Z+\mu\Z)$ sends $P=\tau$
to $(4\mu-1)/10$. Moreover, we have $A = (\begin{smallmatrix}
-9 & 4 \\ 20 & -9 \end{smallmatrix})$ on the geodesic $\{\frac25, \frac12\}$,
and $A = (\begin{smallmatrix} 11 & -4 \\ -30 & 11 \end{smallmatrix})$
on $\{\frac25, \frac13\}$. By Lemma \ref{lem:cc2}, it follows that the
complex conjugation on $\C/(\Z+\mu\Z)$ is given by $z \mapsto -\bar z$,
and a similar computation as above gives
\begin{equation} \label{eq regint 10-2}
\langle \gamma_\infty, S_{P,s} \rangle \sim_{|a| \to \infty} -2\pi \frac{10^3}{3} y_\mu B_3\Bigl(\Bigl\{\frac{2s}{5}\Bigr\}\Bigr),
\end{equation}
where $\gamma_\infty$ denotes a generator of $H_1(E_{t_\infty}, \Z)^-$.

\mybullet $\tau \to \frac3{10}$.
We have $B = (\begin{smallmatrix} 3 & -1 \\ 1 & 0 \end{smallmatrix})$,
$A = (\begin{smallmatrix} -19 & 6 \\ 60 & -19 \end{smallmatrix})$ on the
geodesic $\{\frac3{10}, \frac13\}$, and $A = (\begin{smallmatrix}
11 & -3 \\ -40 & 11 \end{smallmatrix})$ on $\{\frac3{10}, \frac14\}$.
This gives
\begin{equation} \label{eq regint 10-3}
\langle \gamma_1, S_{P,s} \rangle \sim_{|a| \to \infty} -2\pi \frac{10^3}{3} y_\mu B_3\Bigl(\Bigl\{\frac{3s}{10}\Bigr\}\Bigr),
\end{equation}
where $\gamma_1$ denotes a generator of $H_1(E_{t_\infty}, \Z)^-$.

\mybullet $\tau \to \frac15$.
We have $B = (\begin{smallmatrix} 2 & -1 \\ 1 & 0 \end{smallmatrix})$,
$A = (\begin{smallmatrix} -9 & 2 \\ 40 & -9 \end{smallmatrix})$ for the
geodesic $\{\frac15, \frac14\}$, and $A = (\begin{smallmatrix}
1 & 0 \\ -0 & 1 \end{smallmatrix})$ for $\{\frac15, 0\}$.
This gives
\begin{equation} \label{eq regint 10-4}
\langle \gamma_0, S_{P,s} \rangle \sim_{|a| \to \infty} -2\pi \frac{10^3}{3} y_\mu B_3\Bigl(\Bigl\{\frac{s}{5}\Bigr\}\Bigr),
\end{equation}
where $\gamma_0$ denotes a generator of $H_1(E_{t_0}, \Z)^-$.

We find the limit behaviour of $y_\mu$ in each of these cases
using the expansions of the roots of $f_a(X)$ in the proof of
Lemma \ref{lemma:unit}(2), together with the Fourier expansion of
$t \circ \nu'$ when $\tau$ approaches the cusp under consideration.
We find $y_{\mu} \sim \frac{10\log|a|}{2\pi}$ when $\tau \to \frac1{10}, \frac3{10}$
and $y_{\mu} \sim \frac{5\log|a|}{2\pi}$ when $\tau \to \frac1{5}, \frac25$.
Combining this with \eqref{eq regint 10-1}, \eqref{eq regint 10-2},
\eqref{eq regint 10-3} and \eqref{eq regint 10-4}, we have the limit
formula for the regulator of $2S_{P,1}$, $2S_{P,2}$, $2S_{P,3}$ and $2S_{P,4}$.
The linear independence of these elements follows, noting that the
determinant of $\bigl(B_3(\{\frac{st}{10}\})\bigr)_{1 \le s,t \le 4}$ is non-zero.

\smallskip\noindent
The proof of Theorem~\ref{thm:independence} is now complete.

\section{\texorpdfstring{$K_2^T$}{K2T} of families of elliptic curves over cubic fields} \label{K2-cubic}

In this section, we construct other families of elliptic curves~$ E $ over cubic fields with three elements in~$ \INTE $.
As mentioned in the introduction, the results in this section
are mostly independent of the rest of the paper, being based
directly on Lemma~\ref{lemma:unit}(1) and the regulator pairing~\eqref{regulator-pairing}.

\begin{prop} \label{prop:curvecubic}
Let $F$ be a field of characteristic other than~2,
with~$ p \ne 0, -1 $ in~$ F $.
Let $A=p^2+p+1 $, $B=p^2(p+1)^2$, ~$\l=1,2,3$ or~$4$, and~$ t = \frac{4 B}{\l  A^3} $.
If~$ 27 t \ne 4 $,
then, with~$ h(x) = A x + B $, the normalisation of the curve defined by
\begin{equation} \label{eqn:curvecubic}
y^2 + \bigl ( 2 x^3 + \l h(x)^2 \bigr) y + x^6 = 0
\end{equation}
is an elliptic curve~$E$.
With $ C_1 = 6 $, $ C_2 = 4 $, $ C_3 = 3 $, and $ C_4 = 2 $,
the four elements
\begin{equation*}
M = \left\{ - \frac{y} {x^3} , \frac{ h(x) } { h(0) } \right\}
\quad\quad
M_q = C_\l  \left\{ - \frac{y}{x^3}, q^{-2} x+1 \right\}
\hskip-7pt\quad (q = p,  p+1, p(p+1) )
\end{equation*}
are in~$ K_2^T(E)$ and satisfy~$ 2 C_\l M = \sum_q M_q $.
If $ p^2 + p + 1 = 0 $ then $ M = 0 $.
\end{prop}

\begin{proof}
For any $ \l \ne 0 $, completing the square and letting $ \tilde h(x) = x^3 + \frac {\l} 2 h(x)^2 $,
$x' = \frac{Ax}{B}$ and~$y'=\frac{2(y+\tilde h(x))}{\l Bh(x)} $ leads to an equation
\begin{equation}\label{wei-model}
y'^2 = t x'^3 + (x' + 1)^2
.
\end{equation}
The conditions on~$ p $ and~$ t $ in the proposition are such that the right-hand side is
a cubic without multiple roots, so this equation defines an elliptic curve.

However, in order to check that~$ M $ and the~$ M_q $ are in $ K_2^T(E) $ it is more convenient to work
with~\eqref{eqn:curvecubic},
which we view as defining a cover~$ C $ of $ \P_F^1 $ of degree~2
by combining it with a part given by
\begin{equation} \label{coverpart}
(\tilde y + 1 )^2 + \l  \tilde x ( B \tilde x + A )^2 \tilde y  = 0
\end{equation}
with~$ x = 1/\tilde x $ and $ y = \tilde y / \tilde x^3 $.
This adds one point~$ (\tilde x, \tilde y) = (0,-1) $ at infinity.
We view $ E $ as the normalisation of~$ C $.

The divisors of the functions in~$ M $ involve only points
on~$ E $ lying above~$ (x_0,y_0) $ on~$ C $ with $ x_0 = 0 $
or~$ h(x_0) = 0 $, as well as the point at infinity.
We shall show that at such points on~$ C $, at least one of the functions
in~$ M $ is regular with value~1.
The same then holds at the points on~$ E $ above
them, so that~$ M $ is in~$ K_2^T(E) $.

If~$ x_0 = 0 $ then $ h(x)/h(0) $ is regular with value~1 at~$ (x_0,y_0) $.
If $ h(x_0) = 0 $ then from~\eqref{eqn:curvecubic} we find~$ y_0+x_0^3=0 $,
so that~$ -y/x^3 $ is regular with value~1 at~$ (x_0, y_0) $ because~$ x_0 \ne 0 $ by our conditions on~$ p $.
At the point at infinity, the function~$ - y/x^3 = - \tilde y $
is regular with value~1.

We use the same strategy for~$ M_q = \left\{ (-y/x^3)^{C_\l} , q^{-2} x + 1\right\}$, where the relevant
points~$ (x_0,y_0) $ on~$ C $ are those with $ x_0 = 0 $, $ x_0 = -q^2 $,
and the point at infinity. The argument at this last point is
the same as that used for~$ M $, and for~$ x_0 = 0 $ the
function~$ q^{-2} x + 1 $ is regular with value~1 at~$ (x_0,y_0) $.
Now in the remaining points~$ x_0 = - q^2 \ne 0 $, and~$ (-y/x^3)^{C_\l} $ is regular.
From the polynomial identity
\begin{equation} \label{apol-id}
 x^3 + h(x)^2  =  (x+ p^2) (x+ (p+1)^2) (x + p^2 (p+1)^2 )
\end{equation}
together with~\eqref{eqn:curvecubic} we obtain~$ y_0^2 + (2-\l) x_0^3 y_0 + x_0^6 = 0 $.
Hence the function has value~$ (-y_0/x_0^3)^{C_\l} = 1 $
by our choice of~$ \l $ and $ C_\l $.

Note that if $ p^2 + p + 1 = 0 $ then $ h(x) = h(0) $, and~$ M = 0 $.
So it remains to show that $ \sum_q M_q = 2 C_\l M $.
For this, we see from~\eqref{apol-id} and~$ h(0) = \prod_q q $ that
\begin{align*}
\sum_{q}\left\{ - \frac{y}{x^3}, q^{-2} x+1 \right\}-2M &=
\left\{ - \frac{y}{x^3},  \frac{ x^3+h(x)^2 }{h(x)^2}\right\}\\
&= \left\{ - \frac{y}{x^3}, \frac{ y^2+(2-\l)x^3y+x^6 }{(y+x^3)^2}\right\}\\
&=\left\{- \frac{y}{x^3}, \frac{(y^2+(2-\l)x^3y+x^6)}{y^2}\right\}\\
&= \left\{\frac{1}{a}, 1+(\l-2)a+a^2 \right\} \qquad \left(-\frac{x^3}{y}=a\right)
,
\end{align*}
where we used~$(y+x^3)^2= -\l h(x)^2y$
and~$ y^2+(2-\l)x^3y+x^6 = -\l ( x^3+h(x)^2) y $,
and the identity $ \{y/x^3 , y+x^3 \}  = \{ -y/x^3, y \}+  \{-1 , y+x^3 \} $
was obtained using the Steinberg relation~$ 0 = \{y/(y+x^3) , x^3/(y+x^3) \} $.
The last displayed element is~$C_\l$-torsion by our choice of~$ \l $
and $ C_\l $, so that~$ \sum_{q} M_q = 2C_\l M $.
\end{proof}

\begin{remark} \label{remark:identifications}
The substitutions~$p\rightarrow 1/p$ or $p\rightarrow -(p+1)/p$ in~\eqref{eqn:curvecubic}, together with~$(x,y) \mapsto (x/p^4,y/p^{12}) $,
give isomorphic curves for the same fixed $ \l $.
Under these isomorphisms, the elements~$ M $ correspond, as do the~$M_q$ (up
to reindexing). Repeated applications give such identifications
also for~$p\rightarrow -(p+1)$, $p\rightarrow -1/(p+1)$ and $p\rightarrow -p/(p+1)$.
If we let $ u = p+1 $ (cf.~Theorem~\ref{theorem:cubic} below),
then these transformations of $p$ correspond to $u \rightarrow u^{-1} $, $ 1 - u^{\pm1} $ and~$ (1-u^{\pm1})^{\pm1} $,
as in Remark~\ref{remark:isofields}.
\end{remark}

\begin{remark} \label{rk Z2 Z6}
For any $ \l \ne 0 $ in $ F $ such that~$ 27 t \ne 4 $, so that~\eqref{wei-model} defines
an elliptic curve~$ E $, if we let~$\mu = \sqrt{\l^2-4\l}$,
then on $ E $ we have the points~$ P=(0, 1) $ and
$ T_{q}=(-\frac{Aq^2}{B}, \frac{\mu q^3}{\l B}) $ for $ q $ as in Proposition~\ref{prop:curvecubic}.
The divisor
of~$ y' - (x'+1) $ is equal to~$ 3(P) - 3(O) $, so~the point~$ P$ has order~3.
We also have
\begin{equation*}
- \frac{y}{x^3}
=
\frac{x^3}{ y + 2x^3+\l h(x)^2}
=
\frac{ tx'{}^3 }{ \left( y' + x' + 1 \right)^2 }
=
\frac{y' - (x' + 1) }{y' + (x' + 1) }
,
\end{equation*}
so that
\begin{align}
\begin{split}\label{eqn:divMq}
\left(-\frac{y}{x^3}\right) &= 3(P)-3(-P) \\
\left(q^{-2}x+1\right) &= \left(  \frac{B}{Aq^2} x' + 1 \right) = (T_{q}) + (-T_{q}) - 2(O)
.
\end{split}
\end{align}

For $\l=4$, the $T_{q}$ are $ F $-rational and have order~2 because their $y'$-coordinate is zero,
and from~\eqref{apol-id} we see they are all distinct because
$ x^3 + h(x)^2 $ by assumption has no multiple roots.
So the divisors of the functions in the~$M_q$ are supported in (and
generate) a subgroup of $ F $-rational points of~$ E $ isomorphic
to~$\Z/2\Z \times \Z/6\Z$.

For $\l = 1,2,3$, and~$p$ as in Theorem~\ref{theorem:cubic} below, the norm
of the $ x $-coordinate of each~$3T_{q}$ is a rational function
(not a polynomial) in $ a $, so for general $a$ in $\Z$, it is not an integer.
Using \cite[Theorem~7.1]{Si}, we see that~$ T_{q} $ is in general not a torsion point.
Hence the divisor of the function $ q^{-2} x + 1 $ in~$M_q$ is in general not supported in torsion points.
\end{remark}

\begin{remark}
(1)
With notation as in Remark~\ref{rk Z2 Z6}, the element~$ M $
in Proposition~\ref{prop:curvecubic} is equal to~$ S_{P,1} = T_{P,1,2} $
of Section~\ref{blochsection} on~\eqref{wei-model}, with $ P = (0, 1 ) $.
Because~$ ( y' - (x'+1) ) = 3 (P) - 3 (O) $, and~$ 2P = -P = (0, -1 ) $,
we have
\begin{equation*}
T_{P,1,2}
=
\left\{ \frac{y' - (x'+1)} { - 2 } , \frac{y' + (x' + 1)} { 2 } \right\}
=
\left\{  \frac{y' - (x'+1)} {y' + (x'+1)} , \frac{y' + (x'+1)} { 2 } \right\}
.
\end{equation*}
Adding the Steinberg relation
$ 0 = \left\{  \frac{y' - (x'+1)} {y' + (x'+1)} , \frac { 2 (x'+1) } {y' + (x'+1)} \right\} $
we obtain~$ M $ because we saw in Remark~\ref{rk Z2 Z6} that~$ \frac{y' - (x'+1)} {y' + (x'+1)} = - \frac{y}{x^3} $, and~$ \frac{h(x)}{h(0)} = x'+1$.

(2)
Replacing $y$ with $\l^{1/2}h(x)y-x^3$, the family~\eqref{eqn:curvecubic} becomes
\begin{equation*}
y^2 + \l^{1/2}h(x)y-x^3 = 0
,
\end{equation*}
which, up to replacing $ x $ with $ -x $, is a family with $ \Z/3\Z $-torsion studied in \cite{DJZ} (see,
e.g., (6.1) of loc.~cit.).
On this model
\begin{equation*}
M = \left\{ \frac{y^2} {x^3} , \frac{ h(x) } { h(0) } \right\}
\quad\quad
M_q = C_\l  \left\{ \frac{y^2}{x^3}, q^{-2} x+1 \right\}
\hskip-7pt\quad (q = p,  p+1, p(p+1) ).
\end{equation*}
For $\l=4$, via~\eqref{apol-id} the~$M_q$ are related to
the elements of loc.~cit.~(see~(6.5) there).
\end{remark}

We can now give the main result of this section. Note that by Remarks~\ref{remark:isofields}
and~\ref{remark:identifications}, up to isomorphism, in it we may restrict to
$ f_a(X) = X^3 + a X^2 - (a+3) X + 1 $ with~$ a \ge -1 $ for the
cyclic cubic fields, and
$ f_a(X) = X^3 + a X^2 - (a+1) X + 1 $ with~$ a \ge 0 $ for the
non-cyclic ones.

\begin{theorem}\label{theorem:cubic}
For fixed $ \v, \v' $ in $ \{\pm1\} $ in Lemma~\ref{lemma:unit}(1), let
the corresponding polynomial $ f_a(X) $ define $ F = \Q(u) $ with $ u $ a root of~$ f_a(X) $,
so that~$ u $ is an exceptional unit of~$ \OO_F $.
Let~$ p = u - 1 $ and~$ \l = 1 $, 2, 3 or 4.
Then, unless~$ \v = \v' = 1 $ and $ \l \abs{a + \frac32} = \frac92 $,
the normalisation of the curve defined by~\eqref{eqn:curvecubic}
defines an elliptic curve~$ E $ over~$ F $, and~$ M $ as well
as the~$M_q$ in~Proposition~\ref{prop:curvecubic}
are in~$ \INTE $.
Moreover, the Beilinson regulator~$R=R(a)$ of the $ M_q $ satisfies
\begin{equation*}
\lim_{\abs{a}\rightarrow \infty}\frac{R(a)}{\log^3 \abs{a}} = 16 C_\l^3
,
\end{equation*}
with~$ C_\l $ as in Proposition~\ref{prop:curvecubic}.
In particular, for $\abs{a} \gg 0$, $ M $ and the~$ M_q $ generate a subgroup
isomorphic to~$ \Z^3 $, in which the~$ M_q $ generate a subgroup of index~$ 2 C_\l $.
\end{theorem}

\begin{proof}
Writing~$ \l A^3 ( 27 t - 4) = 4 (27 B - \l A^3 ) $ as polynomial
in $ u = p+1 $ and factorising over~$ \Q $, one sees that one
cannot have both $ 27 t = 4 $ and~$ f_a(u)= 0 $ except when $ \v = \v' = 1 $, and
$ a = -6 $ or~3 for $ \l = 1 $, or
$ a = -3 $ or~0 for $ \l = 3 $, which are the stated exceptions. In all other
cases the normalisation of the curve defined by~\eqref{eqn:curvecubic}
is an elliptic curve by Proposition~\ref{prop:curvecubic}.

Our proof that~$ M $ and the $ M_q $ are in~$ \INTE $ is inspired
by, but different from, that of~\cite[Theorem 8.3(1)]{DJZ},
as we do not only consider singular points in a fibre over a maximal
ideal~$ \PP $ of~$ \OO_F $, but all points of that fibre.

Let~$ \cC $ be the naive model of~$ C $ defined by~\eqref{eqn:curvecubic}
and~\eqref{coverpart}, considered as equations over~$ \OO_F $.
Through iterated blowups (which include normalisations) we
obtain a regular~$ \cE $ with a morphism to~$ \cC $ (see \cite[Corollary~8.3.51]{Liu}).
Its generic fibre~$ E $ is the normalisation of~$ C $, so that $ \cE $ is a regular,
flat and proper model of~$ E $ over~$ \OO_F $.

With~$ M = \{ - y/x^3, h(x)/h(0) \} $ in~$ K_2^T(E) $,
we need to show that~$ T_\DD(M) = 1 $ for each irreducible component~$ \DD $
in each fibre~$ \cE_\PP $ of~$ \cE $ over~$ \PP $.
For~$ \DD $ surjecting to an irreducible component of~$ \cC_\PP $ this follows
if the two functions in~$ M $ are generically defined and non-zero
on that irreducible component of~$ \cC_\PP $.
Using~\eqref{eqn:curvecubic} this is easily checked to hold on
all irreducible components of every~$ \cC_\PP$ for the functions $ x $, $ y $, $ h(x) $ and~$ h(0) $,
because~$ h(0) $ is in~$ \OO_F^\times $.

For an irreducible component~$ \DD $ mapping to a closed point
of~$ \cC_\PP $, we shall show that either the two functions in~$ M $
are regular on~$ \cC $ at the point and have non-zero values
there, or at least one of the two functions is regular on~$ \cC $
at the point with value~1.
Considering the valuations along~$ \DD $ of the pullbacks to $ \cE $ of the functions,
and their values on it, one sees that
either condition implies~$ T_\DD(M) = 1 $.
(The first condition was not considered in the proof of~\cite[Theorem 8.3(1)]{DJZ}.)

Both functions in~$ M $ are regular with non-zero values except
perhaps, in fibres~$ \cC_\PP $, points~$ (x_0, y_0) $ with $ x_0 = 0 $ or $ h(x_0) = 0 $,
and the point at infinity.
If~$ h(x_0) = 0 $ then $ x_0 \ne 0 $ by our conditions on~$ p $, so
the reduction of~\eqref{eqn:curvecubic}
modulo~$ \PP $ shows that~$ y_0 / x_0^3 = -1 $, hence the function~$ -y/x^3 $ is regular with
value~1 at~$ (x_0, y_0) $. If $ x_0 = 0 $ then~$ h(x)/h(0) $ is regular
at the point by our conditions on~$ p $, with value~1.
At the point of infinity, where~$ (\tilde x_0, \tilde y_0) = (0, -1) $,
the function~$ -y/x^3 = - \tilde y $ is regular with value~1.

We use the same approach for~$ M _q = \{ (- y/x^3)^{C_\l} , q^{-2} x + 1 \} $.
The functions~$ x $, $ y $, $ x+q $ and~$ q $, with~$ q $ in~$ \OO_F^\times $,
are generically defined and non-zero on any irreducible component of any~$ \cC_\PP $.
The relevant closed points on~$ \cC $ are, in fibres~$ \cC_\PP $,
the point at infinity, and the points~$ (x_0, y_0) $ with $ x_0 = 0 $ or $ x_0 = -q $.
The point at infinity is dealt with as for~$ M $.
If~$ x_0 = 0 $ then the function~$ q^{-2} x + 1 $ is regular with value~1 at
the point because~$ q $ is in~$ \OO_F^\times $.
If $ x_0 = -q $ then from~\eqref{apol-id} we obtain~$ h(x_0)^2 = - x_0^3 $,
so that~$ y_0^2 + (2-\l) x_0^3 y_0 + x_0^6 = 0 $ by~\eqref{eqn:curvecubic}.
It follows as in the proof of Proposition~\ref{prop:curvecubic} that~$ (-y_0/x_0^3)^{C_\l} = 1 $,
hence the function~$ -y/x^3 $ is regular with value~1 at~$ (x_0, y_0) $.

We now prove the limit result for the Beilinson regulator of the~$ M_q $
using the model defined by~\eqref{wei-model}. We write this
as $ y'{}^2 = S(x') $ with~$ S(x') = t x'{}^3 + (x'+1)^2 $, and
let~$ E_p $ be the elliptic curve it defines. Then
from Remark~\ref{rk Z2 Z6} we see that for~$ q = p $, $ p+1 $ and~$ p(p+1) $,
in $ K_2^T(E_p) $ we have
\begin{equation*}
M_q
=
C_\l \Bigl \{ \frac{ y' - (x' + 1) }{ y' + ( x' + 1) } , \frac{B}{q^2 A} x' + 1 \Bigr \}
=
C_\l \Bigl \{ \frac{  t x'^3}{ \bigl( y' +  x' + 1 \bigr )^2 } , \frac{B}{q^2 A} x' + 1 \Bigr \}
.
\end{equation*}

For an embedding~$ \s : F \to \C $ we know from
the proof of Lemma~\ref{lemma:unit}(1) that, as~$ \abs{a} \to \infty $, $ \s(p) $ is real,
and that precisely one of
\begin{equation} \label{eqn:limit-p-a}
\log \abs{\s(p)} \sim - \log \abs{a}, \quad \log\abs{\s(p)+1} \sim-  \log \abs{a}, \quad \log\abs{\s(p)} \sim \log\abs{a}
\end{equation}
holds. We for now suppress~$ \s $
from the notation, and let~$ p $ be a real number with $ \abs{p} \to 0 $
or~$ \infty $, or $ \abs{p+1} \to 0 $.
For~$ p $ close enough
to the limit, we shall construct a generator~$ \g_p $ of $ H_1(E_p(\C) , \Z)^- $
and evaluate the pairing~$ \langle \g_p, M_q \rangle $ in~\eqref{regulator-pairing}.

Note~$ A, B $ and~$ t $ are real and positive.
So for $ r > 0 $ we have~$ S(r) > 0 $, $ S(-1) < 0 $,
and $ S(-r) > 0 $ is equivalent to~$ t <  r^{-3} ( r - 1)^2 $.
For fixed $ r \ne 1 $ this holds if $ \abs{p} $ or
$ \abs{p+1} $ is close enough to 0 or $ \infty $. Now fix~$ r > 1$, and consider
only real~$ p $ such that~$ S(-r) > 0 $.
Then~$ S(x') $ has three real roots~$ \a_1 <  \a_2 < \a_3 $,
which satisfy~$ \a_1 < -r < \a_2 < -1 < \a_3 < r $.
From this one easily sees that~$H_1(E_p(\C), \Z) = H_1(E_p(\C), \Z)^- \oplus H_1(E_p(\C), \Z)^+ $,
and that a loop~$ \g_p $ in~$ E_p(\C) $ lifting the circle in~$ \C $
with~$ \abs{x'} = r $ generates~$ H_1(E_p(\C), \Z)^-$
(cf.~ the discussion preceding Figure~1(a) on \cite[p.361]{DJZ}).
In fact, for~$ \abs{x'} = r > 1 $ we have~$ \abs{ t x'^3/( x' + 1 )^2} < 1 $
because~$ t < r^{-3} (r-1)^2 $, and
rewriting~\eqref{wei-model} as~$ \bigl ( y'/( x' + 1 ) \bigr )^2 = 1 + t x'^3/( x' + 1 )^2 $,
we can use~$ y'/( x' + 1 ) = \bigl( 1 + t x'^3/( x' + 1 )^2 \bigr)^{1/2} = 1 + \frac12 t x'^3/( x' + 1 )^2 + \dots $
to make the $ y' $-coordinate on~$ \g_p $ explicit as function of~$ x' $.

In order to evaluate~$ \langle \g_p, M_q \rangle $,
we now write~$ M_q $ as~$ C_\l \Bigl \{ \frac{ t x'^3 }{ (2x'  + 2)^2 f^2} , \frac{B}{q^2 A} x' + 1 \Bigr \} $
with~$ f =  2^{-1} \bigl( 1+ y'/( x' + 1 )\bigr ) $,
and compute the following integrals.

\mybullet
$ \frac1{2\pi} \int_{\g_p} \eta( t x'^3 (2 x'  + 2)^{-2}  , \frac{B}{q^2 A} x' + 1 ) $.
This can be done in~$ \C $ using~$ \abs{x'} = r $.
With the 1-form closed, it is a residue calculation in~$ \C $,
the residues occurring at~$ x' = 0 $, $ -1 $ and
$ - q^2 A/B $. The first contributes~0, the second always contributes,
and the third contributes only if~$ q^2 A/B < r $.
But~$ q^2 A/B =  q^2 (p^2+p+1)  p^{-2} (p+1)^{-2} $ has the following limits for~$ q = p $, $ p+1 $, $ p(p+1) $:
$ \infty $, 1, 1 as~$ p \to -1 $;
1, $ \infty $, 1 as~$ p \to 0 $; and
1 , 1, $ \infty $ as~$ \abs{p} \to \infty $.

So if this limit is~$ \infty $, then only the residue at~$ x' = -1 $
contributes. The result is~$ 2 \log\abs{1-B/(q^2A)} $,
which approaches~0 because $ B/(q^2 A) \to 0 $.
In the other six cases the limit is~1, and the residue at~$ x' = - q^2 A /B $ also
contributes. Using~$ t = \frac{4 B}{\l  A^3} $ we see
the value of the integral equals~$ 2 \log(q/A) - \log(\l) $.
Because~$ A = p^2 + p + 1 $ we have~$ \log(q/A) - \log\abs{p+1} \to 0 $ in the two cases where $ p \to -1 $,
$ \log(q/A) - \log\abs{p} \to 0 $ in those where~$ p \to 0 $
and
$ \log(q/A) + \log\abs{p} \to 0 $ in those where~$ p \to \infty $.

\mybullet
$ \frac1{2\pi} \int_{\g_p} \eta( f , \frac{B}{q^2 A} x' + 1 ) $.
This approaches~0 when~$ \abs{p} \to 0 $ or~$ \infty $, or~$ \abs{p+1} \to 0 $, as
one sees from the expansion of~$ y'/( x' + 1 ) $, in which~$ t \to 0 $ in
all cases, and the limit behaviour of~$ q^2 A B^{-1} $.

Reinstating~$ \s $ into the notation, we see by
combining the above with~\eqref{eqn:limit-p-a} that
\begin{equation*}
 \lim_{\abs{a}\rightarrow \infty}\frac{\left(\langle \g_{\s(p)} , M_{\s(q)} \rangle\right)_{\s , q} }{ - 2 C_\l \log \abs{a}}  =
\left(
  \begin{array}{ccc}
   0  & 1 & 1 \\
   1  & 0 & 1 \\
   1  & 1 & 0 \\
  \end{array}
\right).
\end{equation*}
Taking the absolute value of the determinants gives the limit behaviour of the
Beilinson regulator of the~$M_q$, which implies
that for~$\abs{a}\gg 0$ the $ M_q $ generate a subgroup isomorphic
to~$ \Z^3 $. That this also holds if we add~$ M $
as a generator, with the subgroup generated by the~$ M_q $ of
index~$ 2 C_\l $, then follows from the relation~$ \sum_q M_q = 2 C_\l M $.
\end{proof}

\begin{remark}
One can also prove the limit statement in Theorem~\ref{theorem:cubic} along
the lines of the proof of~\cite[Theorem~6.6]{LJ15}, using~\cite[Lemma~3.4]{LC18}, which generalises~\cite[Lemma~6.4]{LJ15}.
This involves some rewriting and additional estimates. But
with~$ y' $ quadratic over~$ x' $ in this case,
our current approach results in a much more explicit formula for the regulator pairing,
based mostly on a residue calculation in the complex plane.
This gives a different point of view, 
and using estimates for the only unknown term in the regulator
pairing, it is possible to determine explicitly for which~$ \abs{a} $ sufficiently large
the Beilinson regulator of the~$ M_q $ is non-zero.
\end{remark}

\section{Numerical results} \label{numerics}
Here we use PARI/GP \cite{PARI}
to numerically verify Beilinson's conjecture for elliptic curves,
as formulated in Section~\ref{conjecture}, in some
of the families in Sections~\ref{K2-families} and~\ref{K2-cubic}.

If in~\eqref{eqn:L-reg} we use $ m $ linearly independent elements of $ \INTE $ modulo torsion
instead of its~$ \Z $-basis~$ \a_1,\dots,\a_m $ then~$ Q $ is
replaced with~$ \tQ = Q / l $ for~$ l $ the index of the
subgroup
generated by the elements.
Note that the numerator of such a (computable)~$ \tQ $ divides that
of~$ Q $, even though~$ Q $ itself remains unknown.
We hope that in the future there will be a conjecture for the
$ K $-groups of curves over number fields, in the
spirit of Lichtenbaum's conjecture for the $ K $-groups of number
fields~\cite{lichtenbaum},
and that the examples in the tables (with various large prime
factors in the numerator of some~$ \tQ $, and hence of~$ Q $) can be useful
for formulating or numerically verifying such conjectures (cf.~the beginning of~\cite[Remark~10.14]{DJZ}).

\begin{table}[tbh]
\caption{Data for $N=7$.}\label{table7}
\centering
\begin{tabular}{|c|c|c|c|c|}
\hline \vph
$a$ & $d$ & $c$ & $L^*(E,0)$ & $\tQ$  \\
\hline \vph
$-15$ & $47 \cdot 911$ & $1259$ & $-33502189.0313992549$ & $-2^{2} \cdot 3 \cdot 5 \cdot 7^{-5} \cdot 11 \cdot 19^{2}$\\
$-14$ & $32009$ & $2^{3} \cdot 113$ & $-10623031.4445936662$ & $-2^{4} \cdot 3 \cdot 7^{-4} \cdot 13 \cdot 19$\\
$-13$ & $97 \cdot 241$ & $617$ & $5120119.76612789544$ & $2^{4} \cdot 3^{3} \cdot 7^{-5} \cdot 103$\\
$-12$ & $17 \cdot 977$ & $2^{3} \cdot 7^{2}$ & $-1880894.82849877177$ & $-2^{4} \cdot 3 \cdot 5 \cdot 7^{-4} \cdot 11$\\
$-11$ & $7^{2} \cdot 233$ & $223$ & $-298039.373168375865$ & $-2^{3} \cdot 3^{2} \cdot 7^{-5} \cdot 47$\\
$-10$ & $7537$ & $2^{3} \cdot 13$ & $97160.8374079165014$ & $2^{2} \cdot 7^{-4} \cdot 47$\\
$-9$ & $4729$ & $29$ & $9629.57910342935362$ & $2^{2} \cdot 7^{-5} \cdot 41$\\
$-8$ & $2777$ & $2^{3}$ & $654.969737972166692$ & $2^{3} \cdot 7^{-5}$\\
$-7$ & $1489$ & $13$ & $-434.612790919979199$ & $-2^{3} \cdot 7^{-5}$\\
$-6$ & $17 \cdot 41$ & $2^{3}$ & $62.5967792926624719$ & $7^{-5}$\\
$-5$ & $257$ & $7^{2}$ & $32.3160033216808382$ & $7^{-5}$\\
$-4$ & $7^{2}$ & $2^{3} \cdot 13$ & $4.00033992864846188$ & $3^{-1} \cdot 7^{-5}$\\
$-3$ & $-23$ & $167$ & $1.28519307117788975$ & $2^{-1} \cdot 7^{-6}$\\
$-2$ & $-31$ & $2^{3} \cdot 29$ & $2.68728590621394687$ & $7^{-6}$\\
$-1$ & $-23$ & $293$ & $2.34067746029050612$ & $7^{-6}$\\
$0$ & $-23$ & $2^{3} \cdot 43$ & $2.81741246119387006$ & $7^{-6}$\\
$1$ & $-31$ & $379$ & $-3.65888851182820071$ & $-7^{-6}$\\
$2$ & $-23$ & $2^{3} \cdot 7^{2}$ & $3.20759739648506351$ & $7^{-6}$\\
$3$ & $7^{2}$ & $13 \cdot 29$ & $14.5301315201187081$ & $7^{-5}$\\
$4$ & $257$ & $2^{3} \cdot 41$ & $235.760168840014734$ & $7^{-4}$\\
$5$ & $17 \cdot 41$ & $239$ & $1671.96067772426875$ & $2 \cdot 3 \cdot 5 \cdot 7^{-5}$\\
$6$ & $1489$ & $2^{3} \cdot 13$ & $4051.92834496448134$ & $7^{-3}$\\
$7$ & $2777$ & $83$ & $-6590.94375552556550$ & $-2 \cdot 5 \cdot 7^{-5} \cdot 11$\\
$8$ & $4729$ & $2^{3} \cdot 41$ & $114693.828270615380$ & $2^{3} \cdot 3^{3} \cdot 7^{-4}$\\
$9$ & $7537$ & $7^{2} \cdot 13$ & $520366.913326434323$ & $2 \cdot 3 \cdot 7^{-4} \cdot 137$\\
$10$ & $7^{2} \cdot 233$ & $2^{3} \cdot 127$ & $-1485239.71027494934$ & $-2 \cdot 3^{2} \cdot 7^{-4} \cdot 113$\\
$11$ & $17 \cdot 977$ & $1471$ & $5790649.98684165696$ & $2^{4} \cdot 3 \cdot 5^{2} \cdot 7^{-5} \cdot 41$\\
$12$ & $97 \cdot 241$ & $2^{3} \cdot 251$ & $17255203.9121322960$ & $2^{4} \cdot 3^{2} \cdot 7^{-4} \cdot 131$\\
$13$ & $32009$ & $2633$ & $28504752.7830982117$ & $2^{8} \cdot 3 \cdot 7^{-4} \cdot 37$\\
$14$ & $47 \cdot 911$ & $2^{3} \cdot 419$ & $93361926.2369695039$ & $2^{3} \cdot 3 \cdot 7^{-4} \cdot 3571$\\
$15$ & $73 \cdot 769$ & $43 \cdot 97$ & $192572866.057081271$ & $2^{3} \cdot 3^{2} \cdot 7^{-4} \cdot 43 \cdot 53$\\
\hline
\end{tabular}
\end{table}

\begin{table}[tbh]
\caption{Data for $N=8$.}\label{table8}
\centering
\begin{tabular}{|c|c|c|c|c|}
\hline \vph
$a$ & $d$ & $c$ & $L^*(E,0)$ & $\tQ$  \\
\hline \vph
$0$ & $-23$ & $7 \cdot 23$ & $1.22348267827696836$ & $2^{-23}$\\
$1$ & $-31$ & $3^{3} \cdot 17$ & $6.64291810562928558$ & $2^{-21}$\\
$2$ & $-23$ & $5 \cdot 137$ & $5.97110504152047155$ & $2^{-21}$\\
$3$ & $7^{2}$ & $7 \cdot 113$ & $31.2948786232840397$ & $2^{-18}$\\
$4$ & $257$ & $3^{3}$ & $25.2202129687784361$ & $2^{-18} \cdot 3^{-1}$\\
$5$ & $17 \cdot 41$ & $11 \cdot 41$ & $3130.70411060858445$ & $2^{-15} \cdot 3$\\
$6$ & $1489$ & $7 \cdot 13$ & $3377.15438740388289$ & $2^{-13}$\\
$7$ & $2777$ & $3^{3} \cdot 5 \cdot 7$ & $-110191.314028644712$ & $-2^{-10} \cdot 3$\\
$8$ & $4729$ & $17 \cdot 127$ & $806249.659144856084$ & $2^{-13} \cdot 11 \cdot 13$\\
$9$ & $7537$ & $19 \cdot 199$ & $-3399020.63508445448$ & $-2^{-12} \cdot 257$\\
$10$ & $7^{2} \cdot 233$ & $3^{3} \cdot 7 \cdot 31$ & $9860642.47040826474$ & $2^{-11} \cdot 3 \cdot 109$\\
$11$ & $17 \cdot 977$ & $23 \cdot 367$ & $-38313626.2137679483$ & $-2^{-13} \cdot 4547$\\
$12$ & $97 \cdot 241$ & $5 \cdot 463$ & $22214626.7118122391$ & $2^{-14} \cdot 4787$\\
$13$ & $32009$ & $3^{3} \cdot 7$ & $2759510.81590883242$ & $2^{-13} \cdot 3 \cdot 7 \cdot 13$\\
$14$ & $47 \cdot 911$ & $7 \cdot 29 \cdot 97$ & $-549654076.156923184$ & $-2^{-12} \cdot 3^{4} \cdot 311$\\
$15$ & $73 \cdot 769$ & $17 \cdot 31 \cdot 47$ & $1205314746.12464172$ & $2^{-9} \cdot 5 \cdot 1289$\\
\hline
\end{tabular}
\end{table}

\begin{table}[htb]
\caption{Data for $N=10$.}\label{table10}
\centering
\begin{tabular}{|c|c|c|c|c|}
\hline \vph
$a$ & $d$ & $c$ & $L^*(E,0)$ & $\tQ$  \\
\hline \vph
$-7$ & $2^{3} \cdot 41^{2}$ & $2^{2} \cdot 23^{2}$ & $67284.5712909244205$ & $2^{-11} \cdot 5^{-5}$\\
$-6$ & $2^{6} \cdot 7^{2} \cdot 37$ & $3^{4} \cdot 7^{2}$ & $12809909.2599370080$ & $2^{-9} \cdot 5^{-4} \cdot 13$\\
$-5$ & $2^{3} \cdot 13 \cdot 17^{2}$ & $2^{2} \cdot 19^{2}$ & $321613.252539691824$ & $2^{-10} \cdot 5^{-4}$\\
$-4$ & $2^{8} \cdot 17$ & $17^{2}$ & $1308.96784301967823$ & $2^{-10} \cdot 5^{-7}$\\
$-2$ & $2^{6} \cdot 5$ & $13^{2}$ & $3.90265959107592883$ & $2^{-14} \cdot 5^{-9}$\\
$-1$ & $2^{3} \cdot 7^{2}$ & $2^{2} \cdot 11^{2}$ & $18.1524378610645748$ & $2^{-14} \cdot 5^{-8}$\\
$0$ & $2^{8}$ & $3^{4}$ & $1.29080207928400602$ & $2^{-14} \cdot 3^{-2} \cdot 5^{-8}$\\
$1$ & $2^{3} \cdot 7^{2}$ & $2^{2} \cdot 7^{2}$ & $7.41655915683319223$ & $2^{-15} \cdot 5^{-8}$\\
$2$ & $2^{6} \cdot 5$ & $5^{2}$ & $0.604505751430063810$ & $2^{-14} \cdot 5^{-10}$\\
$4$ & $2^{8} \cdot 17$ & $7^{2}$ & $211.227406732423650$ & $2^{-11} \cdot 5^{-7}$\\
$5$ & $2^{3} \cdot 13 \cdot 17^{2}$ & $2^{2}$ & $825.817965343090665$ & $2^{-11} \cdot 5^{-7}$\\
$6$ & $2^{6} \cdot 7^{2} \cdot 37$ & $3^{4}$ & $272030.854985666477$ & $2^{-9} \cdot 3^{2} \cdot 5^{-6}$\\
$7$ & $2^{3} \cdot 41^{2}$ & $2^{2} \cdot 5^{4}$ & $111421.646021166774$ & $2^{-10} \cdot 5^{-5}$\\
\hline
\end{tabular}
\end{table}

In practice, we calculate a regulator determinant as in~\eqref{eqn:L-reg} based on~\eqref{eqn:reg}
for elliptic curves over the complex numbers.
The latter involves
the Fourier expansion~\eqref{eqn:Fourier} of the functions~$D_\t$ and~$J_\t$,
which converge exponentially fast,
so that we can compute this determinant numerically to high precision.

Now let $ E $ be an elliptic curve over a number field~$ F $,
and fix~$ \a $ in $ K_2^T(E) $. Let~$ X = \coprod_\s E^\s(\C) $
where $ \s $ runs through the embeddings of~$ F $ into~$ \C $.

For a real embedding $\sigma$, let $ \Lambda $ be the period lattice of a non-zero 1-form on $ E^{\sigma} $, so
$ \overline{\Lambda} = \Lambda $ and the map $ E^{\sigma}(\C) \cong \C / \Lambda $
is compatible with complex conjugation. If we scale $ \Lambda $
using a non-zero real number to $ \Z +\Z \t $, then we may assume
$ \t = x_\t + i y_\t $ is in~$ i \R_{>0} $ or $ \t $ is in $ \frac12 + i \R_{>0} $.
Then under the identification $ E^{\sigma}(\C) \cong \C / \Lambda \cong \C/ (\Z + \Z \t ) $,
$ H_1(E^\sigma(\C), \Z)^- $ is generated by $ \g = [0, \t] $
in the
first case and by~$ \g = [0, 2i y_\tau] $ in the second.
Hence $ \frac{\Omega_\g}{y_\t \Omega_0}= ci $  with $ c $ equal to 1 or 2 in these cases.
Let $u= u_\s$ be the diamond operator applied to $\a^\s$.
By~\eqref{eqn:reg}, we have
\begin{equation} \label{eqn:reg-pairing-diamond}
\langle \g,\a \rangle_X = \langle \g,\a^\s \rangle = -\frac{1}{2\pi }\im(c i R_\t(u)) = -\frac{c}{2\pi}D_\t(u).
\end{equation}

For a pair of conjugated complex embeddings $\s$ and $\bar{\s}$, let $ E^{\s}(\C) \cong \C/ (\Z + \Z \t ) $ and $ E^{\bar{\s}}(\C) \cong \C/ (\Z + \Z \bar{\t} ) $ be compatible with complex conjugation. Then~$\g_1=[0,1]$ and $\g_2=[0,\t]$ form a $\Z$-basis of $H_1(E^\s(\C),\Z)$,
and~$\g_1^-=\g_1 - \overline{\g_1}$ and $\g_2^-=\g_2 - \overline{\g_2}$ become part of a basis of $ H_1(X, \Z)^- $.
Then from the behaviour of the argument under complex conjugation,
we see from~\eqref{regulator-pairing} that $ \langle \g_j^-, \a \rangle_X = 2 \langle \g_j, \a^\s \rangle$,
with the latter computed on $ E^\s(\C) $.
So if $u$ is the result of applying the diamond operator~$ \diamond $ to~$\a^\s$,
then by~\eqref{eqn:reg} we have
\begin{align*}
\langle \g_1^-,\a \rangle_X & = \frac{1}{\pi y_\t}J_\t(u), \\
\langle \g_2^-,\a \rangle_X & = -\frac{1}{\pi y_\t}\im(\t R_\t(u)) = -\frac{1}{\pi}D_\t(u) + \frac{x_\t}{\pi y_\t}J_\t(u)
.
\end{align*}
Since we calculate the determinant, we may ignore the term~$ \frac{x_\t}{\pi y_\t}J_\t(u)$ in $\langle \g_2^-,\a \rangle_X$.
We can find~$ \t $ with $\im(\t)>0$ numerically as~$ \omega_1/\omega_2$, with~$\Z\omega_1+\Z\omega_2$
the period lattice as output by~\texttt{ellperiods} in~PARI/GP. When $E$ is defined over $\R$, $\re(\t)$ equals~0 or $\frac{1}{2}$.
Also, for~$ E $ defined over~$ \C $, given a point on~$ E $,
the function $\texttt{ellpointtoz}$ provides an element of~$ \C $ that corresponds to the point in $ \C $ modulo
the lattice.

We use the function $\texttt{lfun}$ to calculate $ L^*(E,0) $,
with~18 significant digits. This usually requires a
lot of memory when the conductor norm of the elliptic curve is large,
which in practice limits how many examples we can give in each
table. We use the function \texttt{lindep} to recognise the resulting~$ \tQ $
as a non-zero rational number.

\begin{example} \label{ex:N}
For~$N=7$, 8 and~10, let us consider the elliptic curves $E$ over a number field~$F$ with a
point of order~$N$, as described in Theorem~\ref{thm:independence}.

For $N=7$ and~8, let $F$ be defined by $f_a(X)$ with~$ \v = 1 $ and $ \v' = -1 $
as in Lemma~\ref{lemma:unit}(1). So~$u$ is a root of~$ f_a(X) =X^3+aX^2-(a+1)X+1 $
with~$a$ in~$\Z$. By the
lemma, $F$ is non-Abelian except for~$a=-4$ or~$3$, when $F$ is cyclic,
and from its discriminant as computed in the proof one sees~$ F $
is totally real unless~$ a = -3, \dots, 2 $.

In Theorem \ref{thm:independence}, for $ N = 8 $ we have $t=1/(u+1)$. If $a\rightarrow -(a+1)$, then $u \rightarrow 1/u$ and $t \rightarrow 1-t$,
so by Remark~\ref{remark:isoN} we only need to consider $a\geqslant 0$ here.

For $N=10$, let $F$ be defined by $f_a(X)=X^4+aX^3-aX+1$ with $a$ in~$ \Z$,
which is one of the families in Lemma~\ref{lemma:unit}(2)
(see Table~\ref{table:quarticfields}).
By writing out the 28 exceptions mentioned there, we
see that~$ f_a(X) $ is irreducible if and only if $|a| \neq 3$.
By making explicit the end of the proof of Lemma~\ref{lemma:unit}(2), we see the Galois group for
such irreducible $ f_a(X) $ is $D_4$ for $a \neq 0$,  while for $a=0$ it is $\Z/2\Z \times \Z/2\Z$.
Using $f_a(X)=(X^2+\frac a2 X-1)^2-(\frac{a^2}4 - 2)X^2$, it is straightforward to see
that $F$ is totally real unless~$ a = -2, \dots, 2 $, where it has two complex places.
By Theorem~\ref{thm:independence} we know that~$ S_{P,s} $ is in~$ \INTE $ for~$ N = 7 $,
and~$ 2 S_{P,s} $ is for~$ N = 8$ or~10.
We calculate~\eqref{eqn:reg-pairing-diamond} for each embedding using~\eqref{Spdiamond}.
In Tables~\ref{table7} (resp. \ref{table8} and \ref{table10}), for some small $a$ we list $\tQ$ as in \eqref{eqn:L-reg} using~$S_{P,s}$ for~$N=7$ (resp. $2S_{P,s}$ for $N=8, 10$)
with~$1\leqslant s\leqslant \lfloor \frac{N-1}{2}\rfloor$, together with the discriminant $d$ of $F$, conductor norm $c$ of $E$ and $L^*(E,0)$.
\end{example}

\begin{remark} \label{remark:ST}
The simple shape of~\eqref{Spdiamond} is convenient
in the part of the proof of Theorem~\ref{thm:independence}
in Section~\ref{section:independence}. But as used in the proofs of
Theorems~\ref{thm:integrality} and~\ref{thm:independence},
the~$ S_{P,s} $ lie in the subgroup of~$ \INTE $ generated by the~$ T_{P,s,t} $ for
$ N = 7 $, and the~$ 2 S_{P,s} $ in that generated by the~$ 2 T_{P,s,t} $ for $ N = 8 $ and~10.
We can use those, and other elements, including those in~Remark~\ref{remark:fafb},
to obtain Beilinson regulators for which the corresponding~$ \tQ $
is a positive integer multiple of that of the~$ S_{P,s} $ or~$ 2 S_{P,s} $.

Suppose for~$ N = 7 $ that~$ S_{P,1} $, $ S_{P,2} $ and~$ S_{P,3} $ have Beilinson regulator~~$ R_S \ne 0 $.
The~$ S_{P,s} $ lie in the subgroup of~$ \INTE $ generated by the~$ T_{P,s,t} $,
which contains elements (e.g., $ T_{P,1,2} $, $ T_{P,1,3} $
and~$ T_{P,1,4} $)
with Beilinson regulator~$ 7^{-2} R_S $.
(This can be computed using their images
under the diamond operator in~$ \Z[E(\C)]^- $
modulo torsion, which all lie in the rank 3 subgroup generated
by the classes of~$ (P) $, $ (2P) $ and $ (3P) $.)
Using instead the~$ \b_{P,b,c} $ of Remark~\ref{remark:fafb}
(e.g., $ \b_{P, 2, 3} $, $ \b_{P, 2, 4} $ and~$ \b_{P, 2, 5} $)
one obtains a Beilinson regulator that is~$ 7^{-4} R_S $.
But the obstruction mentioned at the end of Remark~\ref{remark:fafb}
is in~$ \OO_F^\times / (\OO_F^\times)^7 $, which is isomorphic
to~$ \Z/7\Z $ for~$ a = -3, \dots, 2 $, and to~$ (\Z/7\Z)^2 $
otherwise. Using this we can
enlarge the group, and obtain regulators~$ 7^{-6} R_S $ for~$ a = -3, \dots,  2 $,
and~$ 7^{-5} R_S $ otherwise.
The resulting~$ \tQ $ are then~$ 7^6 $ or~$ 7^5 $ times that in Table~\ref{table7}.
The limit statement in Theorem~\ref{thm:independence} can
be similarly modified by changing the elements.

Similarly, for $ N = 8 $, we can supplement the~$ 2 S_{P,s} $
with all $ 2 T_{P, s, t} $, and all~$ T_{P,s,t} $ with both~$ s $ and~$ t $
even.
And with $ 2P $ of order~4 and hitting the 0-component in every fibre
of the minimal regular model of the curve, we can also use~$ S_{2P,1} $
as well as~$ T_{2P,1,2} $, $ T_{2P, 1, 3} $ and~$ T_{2P,2,3} $.
In the resulting subgroup of~$ \INTE $ one can find three elements
for which~$ \tQ $ is multiplied by~$ 2^{10} $ compared to that
for the $ 2 S_{P,s} $.

For $ N = 10 $, one can similarly supplement the~$ 2 S_{P,s} $
with various elements: the~$ 2 T_{P,s,t} $, and~$ T_{P,s,t} $
if both~$ s $ and~$ t $ are even,~$ S_{2P,1} $ and~$ S_{2P,2} $
(with $ 2P $ of order~5 and hitting the 0-component at all primes),
$ T_{2P, s, t} $ with~$ 1 < s < t < 5 $,
as well as the~$ \b_{2P,a,b} $ with $ 1 < a < b < 5 $ of
Remark~\ref{remark:fafb}. In the resulting
subgroup of~$ \INTE $ one can find four elements for which~$ \tQ $ is~$ 2^{10} 5^4 $
times that for the $ 2 S_{P,s} $.
\end{remark}

\begin{remark}
As observed in Remark~\ref{remark:isofields}, the isomorphism~$ \Q(u_a) \simeq \Q(u_{-a-1}) $ of cubic fields in Example~\ref{ex:N},
which is used in the example when~$ N = 8 $, does
not necessarily lead to isomorphic curves in the example when~$ N = 7 $; in
fact, in Table~\ref{table7} the values of~$ L^*(E,0) $ do not match.
The same applies in the example when~$ N = 10 $ for the isomorphisms of the quartic
fields with $ a $ and~$ -a $ given by~$ u_a \mapsto -u_{-a} $
or~$ u_{-a}^{-1} $, as can be seen from Table~\ref{table10}.
(The isomorphism for~$ N=10 $ in Remark~\ref{remark:isoN} corresponds
to the automorphism of the field given by~$ u_a \mapsto - u_a^{-1} $.)
\end{remark}

\begin{table}[tb]
\caption{Data for $\l=1$.}\label{tab:l1}
\centering
\begin{tabular}{|c|c|c|c|c|}
\hline \vph
$a$ & $d$ & $c$ & $L^*(E,0)$ & $\tQ$  \\
\hline \vph
$0$ & $-23$ & $2^{3} \cdot 17 \cdot 107$ & $132.724179260406391$ & $2^{-4}\cdot 3$\\
$1$ & $-31$ & $2^{3} \cdot 3^{4} \cdot 17$ & $168.814511547175067$ & $2^{-4}\cdot 3$\\
$2$ & $-23$ & $2^{3} \cdot 19 \cdot 37$ & $53.4019469956784239$ & $2^{-4}$\\
$3$ & $7^{2}$ & $2^{3} \cdot 127$ & $37.1776384769406512$ & $2^{-4} \cdot 3 \cdot 7^{-1}$\\
$4$ & $257$ & $2^{3} \cdot 3^{4}$ & $-721.242054102691853$ & $-2^{-3} \cdot 3$\\
$5$ & $17 \cdot 41$ & $2^{3} \cdot 19$ & $1414.02549043158906$ & $2^{-2} \cdot 3$\\
$6$ & $1489$ & $2^{3} \cdot 17 \cdot 19$ & $83163.7726064265207$ & $2 \cdot 3^{3}$\\
$7$ & $2777$ & $2^{3} \cdot 3^{4} \cdot 37$ & $2915249.85675393311$ & $2^{2} \cdot 3^{3} \cdot 13$\\
$8$ & $4729$ & $2^{3} \cdot 71 \cdot 163$ & $33679082.6389894579$ & $2 \cdot 3^{3} \cdot 241$\\
$9$ & $7537$ & $2^{3} \cdot 37 \cdot 863$ & $260954243.280987485$ & $2 \cdot 3^{3} \cdot 1567$\\
\hline
\end{tabular}
\end{table}

\begin{table}[tbh]
\caption{Data for $\l=2$.}\label{tab:l2}
\centering
\begin{tabular}{|c|c|c|c|c|}
\hline \vph
$a$ & $d$ & $c$ & $L^*(E,0)$ & $\tQ$  \\
\hline \vph
$0$ & $-23$ & $2^{6} \cdot 11 \cdot 23 \cdot 37$ & $4486.81605627777558$ & $2^{-1} \cdot 3 \cdot 5$\\
$1$ & $-31$ & $2^{6} \cdot 3^{3} \cdot 11 \cdot 13$ & $3599.55769844723823$ & $2^{-1}\cdot 3^{2}$\\
$2$ & $-23$ & $2^{6} \cdot 5^{2} \cdot 59$ & $837.555573566513198$ & $2$\\
$3$ & $7^{2}$ & $2^{6} \cdot 13 \cdot 83$ & $-2498.99534192761051$ & $-3$\\
$4$ & $257$ & $2^{6} \cdot 3^{3} \cdot 37$ & $-64543.3050825583931$ & $-2^{2}\cdot 3^3$\\
$5$ & $17 \cdot 41$ & $2^{6} \cdot 11^{2} \cdot 13 \cdot 23$ & $-16392164.6852019715$ & $-2^{2} \cdot 3^4 \cdot 53$\\
$6$ & $1489$ & $2^{6} \cdot 23 \cdot 47 \cdot 179$ & $437520185.347094640$ & $2^{5}\cdot 3^{2} \cdot 1187$\\
\hline
\end{tabular}
\end{table}

\begin{table}[tbh]
\caption{Data for $\l=3$.}\label{tab:l3}
\centering
\begin{tabular}{|c|c|c|c|c|}
\hline \vph
$a$ & $d$ & $c$ & $L^*(E,0)$ & $\tQ$  \\
\hline \vph
$0$ & $-23$ & $2^{3} \cdot 3^{9} \cdot 19$ & $25300.9847248343307$ & $3 \cdot 17$\\
$1$ & $-31$ & $2^{3} \cdot 3^{11}$ & $-21806.9954627600874$ & $-2^{2} \cdot 3 \cdot 5$\\
$2$ & $-23$ & $2^{3} \cdot 3^{9} \cdot 17$ & $-21113.3123276958079$ & $-2^{2} \cdot 3 \cdot 5$\\
$3$ & $7^{2}$ & $2^{3} \cdot 3^{9}$ & $5601.39536780219401$ & $2^{2} \cdot 3$\\
$4$ & $257$ & $2^{3} \cdot 3^{11} \cdot 19$ & $-26042785.9143510709$ & $-2^{3} \cdot 3^{3} \cdot 233$\\
\hline
\end{tabular}
\end{table}

\begin{table}[tbh]
\caption{Data for $\l=4$.}\label{tab:l4}
\centering
\begin{tabular}{|c|c|c|c|c|}
\hline \vph
$a$ & $d$ & $c$ & $L^*(E,0)$ & $\tQ$  \\
\hline \vph
$0$ & $-23$ & $2^{6} \cdot 5 \cdot 7$ & $19.1718016489393019$ & $2^{-3}$\\
$1$ & $-31$ & $2^{6} \cdot 3^{2}$ & $8.95758063575193728$ & $2^{-3} \cdot 3^{-1}$\\
$2$ & $-23$ & $2^{6} \cdot 5 \cdot 11$ & $-25.4138019939166741$ & $-2^{-3}$\\
$3$ & $7^{2}$ & $2^{6} \cdot 7 \cdot 13$ & $241.273298483854998$ & $2^{-1} \cdot 3$\\
$4$ & $257$ & $2^{6} \cdot 3^{2} \cdot 5$ & $-2647.23969149488937$ & $-3^{2}$\\
$5$ & $17 \cdot 41$ & $2^{6} \cdot 5 \cdot 11 \cdot 17$ & $441097.703795075666$ & $2^{3} \cdot 3^3 \cdot 5$\\
$6$ & $1489$ & $2^{6} \cdot 7 \cdot 13 \cdot 19$ & $-4149007.28165801473$ & $-2^{7} \cdot 3^{2} \cdot 7$\\
$7$ & $2777$ & $2^{6} \cdot 3^{2} \cdot 5 \cdot 7$ & $2423760.93043136419$ & $2 \cdot 3^3 \cdot 73$\\
$8$ & $4729$ & $2^{6} \cdot 11 \cdot 17 \cdot 23$ & $99044008.9977606699$ & $2^{7} \cdot 3^{2} \cdot 11^{2}$\\
$9$ & $7537$ & $2^{6} \cdot 5 \cdot 13 \cdot 19$ & $66308672.9214609161$ & $2^5 \cdot 3^{2} \cdot 7 \cdot 41$\\
$10$ & $7^{2} \cdot 233$ & $2^{6} \cdot 3^{2} \cdot 5 \cdot 7$ & $41156246.2610705047$ & $2^{4} \cdot 3^3 \cdot 107$\\
\hline
\end{tabular}
\end{table}

\begin{example}
We now consider the elliptic curves $E$ with elements~$ M $ and
$ M_q $ ($ q = p, p+1, p(p+1)$) in $ \INTE $,
as described in Theorem \ref{theorem:cubic}.
where~$F=\Q(u)$ for~$ u $ a root of~$f_a(X)=X^3+aX^2-(a+1)X+1 $
with~$a$ in~$ \Z$, and $p=u-1$.
As mentioned before Theorem~\ref{theorem:cubic}, we may assume~$a\geqslant 0$.

With notation as in Remark~\ref{rk Z2 Z6}, by \eqref{eqn:divMq}, in $ \Z[E(\C)]^- $ we have
\begin{equation*}
C_\l \left(\frac{y}{x^3}\right)\diamond (q^{-2}x+1) = 6 C_\l \bigl((P+T_q) + (P-T_q) - 2(P)\bigr)
.
\end{equation*}
From this and the relation~$ M_p + M_{p+1} + M_{p(p+1)} = 2 C_\l M  $
we can calculate~$\tQ$ as in~\eqref{eqn:L-reg} for $M_p$, $ M_{p+1} $ and~$ M $.
We list the result in Tables~\ref{tab:l1}, \ref{tab:l2}, \ref{tab:l3} and \ref{tab:l4}, for some small $a$,
together with the discriminant $d$ of $F$, conductor norm $c$ of $E$ and~$L^*(E,0)$.
\end{example}

\end{document}